\documentclass[11pt]{article}

\usepackage{color}
\usepackage{latexsym}
\usepackage{amssymb}
\usepackage{amsmath,amsfonts,amssymb,theorem,euscript,array,enumerate,amsfonts,mathrsfs}
\usepackage{graphicx}

\newtheorem{Theorem}{Theorem}[part]
\newtheorem{Definition}{Definition}[part]
\newtheorem{Proposition}{Proposition}[part]

\newtheorem{Lemma}{Lemma}[part]
\newtheorem{Corollary}{Corollary}[part]
\newtheorem{Remark}{Remark}[part]

\def\esssup_#1{\underset{#1}{\mathrm{ess\,sup\, }}}
\def\essinf_#1{\underset{#1}{\mathrm{ess\,inf\, }}}

\def \trans{^{\scriptscriptstyle{\intercal}}}

\def \trans{^{\scriptscriptstyle{\intercal }}}

\def \I{\mathbb{I}}
\def \N{\mathbb{N}}
\def \R{\mathbb{R}}

\def \E{\mathbb{E}}
\def \F{\mathbb{F}}

\def \P{\mathbb{P}}

\def \Ac{{\cal A}}
\def \Bc{{\cal B}}

\def \Ec{{\cal E}}
\def \Fc{{\cal F}}

\def \Lc{{\cal L}}
\def \Pc{{\cal P}}

\def \Oc{{\cal O}}

\def \Uc{{\cal U}}
\def \Vc{{\cal V}}

\def \Vc{{\cal V}}

\def \eps{\varepsilon}

\def \ep{\hbox{ }\hfill$\Box$}

\def\reff#1{{\rm(\ref{#1})}}

\def\beqs{\begin{eqnarray*}}
\def\enqs{\end{eqnarray*}}
\def\beq{\begin{eqnarray}}
\def\enq{\end{eqnarray}}

\allowdisplaybreaks

\begin{document}

\title{Ergodicity of robust switching control and \\ 
nonlinear system of quasi variational inequalities\thanks{E. Bayraktar is supported in part by the National Science Foundation under grant DMS-1613170. H. Pham is  supported in part by FiME (Laboratoire de Finance des March\'es de l'Energie) and
the ``Finance et D\'eveloppement Durable - Approches Quantitatives'' Chair. }}

\author{Erhan BAYRAKTAR\thanks{Department of Mathematics, University of Michigan, 
		\sf erhan@umich.edu} ~~~
		Andrea COSSO\thanks{Dipartimento di Matematica, Politecnico di Milano, \sf andrea.cosso@polimi.it}~~~
               Huy{\^e}n PHAM\thanks{Laboratoire de Probabilit\'es et Mod\`eles Al\'eatoires, CNRS, UMR 7599, Universit{\'e} Paris Diderot, and
               CREST-ENSAE,  \sf pham@math.univ-paris-diderot.fr}
             }

\maketitle

\date{}

\begin{abstract}
We analyze the asymptotic behavior  for a system of fully nonlinear parabolic and elliptic quasi variational inequalities. These equations are related to robust switching control problems introduced in 
\cite{BCP14}. We prove that, as time horizon goes to infinity (resp. discount factor goes to zero) the long run average solution to the parabolic system (resp. the limiting discounted solution to the elliptic system)  is characterized by a solution of a nonlinear   system of ergodic variational inequalities. 
Our results hold under a dissipativity condition and without any non degeneracy assumption on the diffusion term. Our approach uses mainly probabilistic arguments and in particular a dual randomized game representation for the solution to the system of variational inequalities.  
\end{abstract}

\vspace{5mm}


\vspace{5mm}


\newpage

\section{Introduction}

\setcounter{equation}{0} 
\setcounter{Theorem}{0} \setcounter{Proposition}{0}
\setcounter{Corollary}{0} \setcounter{Lemma}{0}
\setcounter{Definition}{0} \setcounter{Remark}{0}

Let us consider  the following system of forward parabolic quasi variational inequalities
\beq
\label{HJBintro}
\begin{cases}
\min\Big\{\dfrac{\partial V}{\partial T}  - \inf_{u\in U}\big[\Lc^{i,u}V  + f(x,i,u)\big], \\
V(T,x,i) - \max_{j\neq i}\big[V(T,x,j) - c(x,i,j)\big]\Big\} \ = \ 0, &\!\!\!\!(T,x,i)\in(0,\infty)\times\R^d\times\I_m, \\
V(0,x,i) \ = \ g(x,i), &\!\!\!\!(x,i)\in\R^d\times\I_m,
\end{cases}
\enq
where $\I_m:=\{1,\ldots,m\}$, with $m\in\N\backslash\{0\}$, and  $\Lc^{i,u}$ is the second-order differential operator
\beqs
\Lc^{i,u} V  &=& b(x,i,u).D_x V  + \frac{1}{2}\text{tr}\big(\sigma\sigma\trans(x,i,u)D_x^2 V \big).
\enqs
Here $U$ is a compact subset of $\R^q$, and the assumptions on the measurable functions $b,\sigma$, $c$, 
$f$ and $g$ will be made precise in the next section. 
Equation \reff{HJBintro} turns out to be related to a certain robust switching control problem studied in \cite{BCP14}: 
\beq
V(T,x,i) & := &  \sup_\alpha \inf_\upsilon \E_{x,i} \Big[ \int_0^T  f(X_t,I_t,\upsilon_t) dt  + g(X_T,I_T) \label{defVT} \\
& & \hspace{2cm} - \; \sum_{n\in\N}  c(X_{\tau_n},I_{\tau_n^-},I_{\tau_n}) 1_{\{\tau_n < T\}} \Big], \nonumber 
\enq
where  
\begin{equation}
\label{SDEintro}
\begin{cases}
X_t \ \!\!= \ x + \int_0^t b(X_r,I_r,\upsilon_r) dr + \int_0^t \sigma(X_r,I_r,\upsilon_r)dW_r, &   \nonumber \\
I_t \ = \ i1_{\{0\leq t<\tau_0\}} + \sum_{n\in\N} \iota_n 1_{\{\tau_n \leq t <\tau_{n+1} \}}. &  \nonumber \\
\end{cases}
\end{equation}
The piecewise constant process $I$ denotes the regime value at time $t$, whose evolution is determined by the controller through the switching control 
$\alpha$ $=$ $(\tau_n,\iota_n)_n$, while the process $\upsilon$, decided by nature, brings the uncertainty within the dynamics of  state process $X$ and the model.  The control process sets in which $\alpha$ and $\upsilon$ run over is a key issue when considering  stochastic differential game type problems as in the formulation \reff{defVT}.  In the robust switching control problem, the switching control $\alpha$  is of feedback form, meaning that  it is chosen by the controller based only on the past and present information coming from the state and regime processes, while the control $\upsilon$ is more generally of open-loop form, since nature is assumed to be aware of all information at disposal.  Precise formulation of robust switching problem is given in Section~\ref{SubS:Robust}. 
Also see \cite{BCP14}, where we proved by developing stochastic Perron's method (which was introduced in \cite{ BS13} to analyze stochastic control problems) further that the value function $V$ in 
\reff{defVT} is the unique viscosity solution to \reff{HJBintro}.
The corresponding elliptic system of quasi variational inequalities for any 
$\beta$  $>$ $0$ is 
\begin{equation}
\label{HJB_ellipticintro}
\min\Big\{\beta V^\beta - \inf_{u\in U}\big[\Lc^{i,u}V^\beta + f(x,i,u)\big],V^\beta(x,i) - \max_{j\neq i}\big[V^\beta(x,j) - c(x,i,j)\big]\Big\} \ = \ 0,
\end{equation}
for any $(x,i)\in\R^d\times\I_m$. Similar to the parabolic case,  this system is related to a robust switching control problem but this time over an infinite horizon  with discount factor $\beta$.

Although it is a classical topic in stochastic control, optimal switching problem, in which a controller implements a discrete set of controls, has attracted a renewed interest and generated important developments in applied and financial mathematics. They occur naturally  in investment problems with fixed transaction costs, pair trading problems or in the \emph{real options} and is a more realistic set-up than assuming that the controller exerts controls of infinite variation or controls that accumulate local time. The literature on this topic  is quite large and we refer e.g. to  the recent papers by  \cite{LyvPha07}, \cite{ElHam09}, \cite{pham09}, \cite{MR2676760},    \cite{HamZha10}, \cite{HuTan10}, \cite{EliKha14} for the analysis of optimal switching problems either by dynamic programming  or backward stochastic differential methods, and to   \cite{DucZer01}, \cite{CarLud10}, \cite{HamJea07}, \cite{sonzha09}, \cite{ngopha} for various applications to  finance and  real options in energy markets.   Most of this literature has focused on the situation where the state coefficients are known, while in practice there is uncertainty about their real value, which motivated us to consider the framework \reff{defVT} in line of recent investigations about robust control problem. On the other hand, an interesting application of the switching systems that we analyze appears to be key in proving the convergence rate of numerical schemes for Hamilton-Jacobi-Bellman equations; see e.g. \cite{MR2336272}.

The chief goal of this paper is to extend the results of  \cite{lions_perthame86} and \cite{menaldi_perthame_robin90}, to the fully non-linear and degenerate case: That is, to
investigate the large time asymptotics of the value function of the robust optimal switching problem $V(T,.)$ as $T$ goes to infinity, which is  closely related to the asymptotic behavior of  $V^\beta$, as $\beta$ goes to zero.  We should mention that the asymptotics for stochastic control and related Hamilton-Jacobi-Bellman (HJB) equations have been studied in other settings by many authors since the seminal papers \cite{benfre92} and \cite{arilio98}, and has received a renewed interest, see the lectures of P.L. Lions (2014-2015) at Coll\`ege de France. Recent papers include for instance  
\cite{ichiharaishii08} for  the long time behavior  of  Hamilton-Jacobi equations in a semi-periodic setting, or  \cite{robxin13}, which considered large time behavior of semi-linear  HJB equations with quadratic nonlinearity in gradients by combining PDE and stochastic analysis methods. We mention \cite{humadric14}, which proved a rate of convergence for the solution to the semi-linear HJB equation towards the ergodic equation under a weak dissipativity condition. We refer also to \cite{BarlesSouganidisI}, \cite{BarlesSouganidisII}, \cite{Hynd}, \cite{ichsheu13}, \cite{Nagai}. The case of fully nonlinear HJB equation is studied recently in \cite{CFP14} by means of backward stochastic differential 
equation (BSDE) representation for nonlinear parabolic, elliptic and ergodic equations.

We introduce three novel features to this problem: 
First, we consider diffusion coefficients $b$, and $\sigma$ depending not only on the state  process $x$ but also on the regime values $i$; second, we incorporate robustness in our model by considering dependence of the diffusion coefficients on the open loop control $\upsilon$. From a PDE point of view, this makes the problem fully nonlinear due to the infimum over  $u$ $\in$ $U$ in \reff{HJBintro}. Third, we do not impose a non-degeneracy condition on the diffusion term.

Our main result is to prove, under natural dissipativity conditions, 
the existence of a constant $\lambda$ (not depending on state and regime values $x,i$) such that
\beq \label{convergodique}
\frac{V(T,x,i)}{T} \ \overset{T\rightarrow\infty}{\longrightarrow} \ \lambda, & &  
\beta V^\beta(x,i) \ \overset{\beta\rightarrow0^+}{\longrightarrow} \ \lambda,
\enq
for all $(x,i)\in\R^d\times\I_m$. Moreover, $\lambda$ is  the solution to the {\it ergodic} system of variational inequalities: 
\begin{equation}
\label{HJBergodicswitching}
\min\Big\{ \lambda  - \inf_{u\in U}\big[\Lc^{i,u} \phi  + f(x,i,u)\big], \phi(x,i) - \max_{j\neq i}\big[\phi(x,j) - c(x,i,j)\big]\Big\} \ = \ 0.
\end{equation}
Here the unknown in the ergodic equation \reff{HJBergodicswitching}  
is the pair $(\lambda,\phi)$ with $\lambda$ real number, and $\phi$ a real-valued function  on $\R^d\times\I_m$. 
We also show that under suitable conditions  $\lambda$ is the value to a robust ergodic control problem (see Remark \ref{remergodic}).
In the proof of the ergodic convergence \reff{convergodique}, a crucial step is the derivation of a uniform (in $\beta$) Lipschitz estimate on $V^\beta$, and so the equicontinuity of the family $(V^\beta)_\beta$.  
The main difficulty in our context, with respect to previous related works, is that  we do not in general have any regularity on the solution 
$V^\beta$ to \reff{HJB_ellipticintro}. We also cannot rely on the analog of the robust switching control representation \reff{defVT} for the infinite horizon 
problem.  Indeed, because of the feedback form on the switching control $\alpha$, which may then depend on the initial state value, it is not clear how 
to get suitable Lipschitz properties of $V^\beta$, see Section~\ref{S:Motivation} for a more detailed discussion.  To overcome this issue, we instead provide a dual probabilistic game representation of $V^\beta$ based on randomization of the controls $\alpha$ and $\upsilon$, following the idea originally developed in \cite{khaetal10} and \cite{KP12} for stochastic control problem, for which we refer also to \cite{CC16}, \cite{CCP15}, \cite{CPH17}. This representation then allows us to derive the needed Lipschitz estimate for $V^\beta$.

The rest of the paper is organized as follows. Section \ref{secformul} sets the assumptions and formulates the main results.  We recall in Section 
\ref{S:Parabolic} the connection between the parabolic system of quasi variational inequalities \reff{HJBintro} and robust switching control problem as 
studied in \cite{BCP14}, and shows by a simple probabilistic control representation argument (which seems to be new to the best of our knowledge) the long run average convergence of the solution to the parabolic system to a solution of  the ergodic system when it exists. Section \ref{S:Elliptic} is the core of the paper, and is devoted to the existence of a solution to the ergodic system as well as the convergence of the limiting discounted solution to the elliptic system. 
The main point is to state a dual representation for $V^\beta$. This is achieved by using  a randomization approach and a 
BSDE representation for the penalized solution $V^{\beta,n}$  to the elliptic system of variational inequalities, which then permits to derive 
a dual representation for $V^{\beta,n}$, and consequently for $V^\beta$, and then to obtain the key estimates for $V^\beta$.  Finally, some technical results are deferred to the Appendix.

\section{Formulation of the problem and main results} \label{secformul} 

\setcounter{equation}{0} 
\setcounter{Theorem}{0} \setcounter{Proposition}{0}
\setcounter{Corollary}{0} \setcounter{Lemma}{0}
\setcounter{Definition}{0} \setcounter{Remark}{0}

\subsection{Notation and assumptions}

Let $U$ be a compact subset of a Euclidean space $\R^q$ and $\I_m:=\{1,\ldots,m\}$, with $m\in\N\backslash\{0\}$. 
We begin imposing the following assumptions on the coefficients of the system of parabolic and elliptic quasi variational inequalities  \eqref{HJBintro} and \eqref{HJB_ellipticintro}.

\vspace{2mm}

{\bf (H1)}
\begin{itemize}
\item[(i)] $b\colon\R^d\times\I_m\times U\rightarrow\R^d$ and $\sigma\colon\R^d\times\I_m\times U\rightarrow\R^{d\times d}$ are continuous and satisfy (we denote by $\|A\|$ $=$ $\sqrt{\text{tr}(AA\trans)}$ the Frobenius norm of any matrix $A$):
\[
|b(x,i,u)-b(x',i,u)| + \|\sigma(x,i,u)-\sigma(x',i,u)\| \ \leq \ L_1 |x-x'|,
\]
$\forall\,x,x'\in\R^d$, $i\in\I_m$, $u\in U$, for some positive constant $L_1$.
\item[(ii)] $f\colon\R^d\times\I_m\times U\rightarrow\R$, $g\colon\R^d\times\I_m\rightarrow\R$, and $c\colon\R^d\times\I_m\times\I_m\rightarrow\R$ are continuous and satisfy
\beqs
|f(x,i,u) - f(x',i,u)| + |c(x,i,j) - c(x',i,j)| &\leq & L_2|x-x'|, \\
|g(x,i)| &\leq & M_2(1+|x|^2),
\enqs
$\forall\,x,x'\in\R^d$, $i,j\in\I_m$, $u\in U$, for some positive constants $L_2$ and $M_2$.
\item[(iii)] $g$ satisfies the inequality
\beqs
g(x,i) &\geq & \max_{j\neq i}\big[g(x,j) - c(x,i,j)\big],
\enqs
$\forall\,x\in\R^d$, $i\in\I_m$.
\item[(iv)] $c$ is nonnegative and the {\it no free loop property} holds: for all $i_1,\ldots,i_k\in\I_m$, with $k\in\N\backslash\{0,1,2\}$,  $i_1=i_k$, and $\text{card}\{i_1,\ldots,i_k\}=k-1$, we have
\beqs
c(x,i_1,i_2) + \cdots + c(x,i_{k-1},i_k) \ > \ 0, \qquad \forall\,x\in\R^d.
\enqs
Moreover, we suppose that $c(x,i,i)=0$, $\forall\,(x,i)\in\R^d\times\I_m$.
\end{itemize}

\noindent We also impose the following \emph{dissipativity condition}.

\vspace{2mm}

{\bf (H2)} \hspace{3mm} For all  $x,x'\in\R^d$, $i\in\I_m$, $u\in U$,
\beq
\label{dissipative}
(x-x').(b(x,i,u)-b(x',i,u)) + \frac12 \|\sigma(x,i,u)-\sigma(x',i,u)\|^2 & \le &  \ -\gamma \,|x-x'|^2,
\enq
\qquad for some constant $\gamma>0$.

\begin{Remark}\label{R:coefficients}
{\rm
(i) Equation \eqref{HJBintro} (and similarly \eqref{HJB_ellipticintro}) turns out to be related to a certain robust switching control problem studied in \cite{BCP14}, as explained in Section \ref{SubS:Robust} below (see also Section \ref{S:Motivation}). For this reason, the set of assumptions {\bf (H1)} is the same as in \cite{BCP14}, but for some additional requirements needed to obtain a backward stochastic differential equation representation (first presented in \cite{KP12} for the case of a classical stochastic control problem) and for the polynomial growth condition on $g$, which in \cite{BCP14} is not necessarily of second degree. This latter assumption plays an important role since it allows to exploit an estimate on the second moment of the state process (reported in Lemma \ref{L:X_diss} below), which follows from the dissipativity condition {\bf (H2)}.

\vspace{2mm}

\noindent(ii) As an example of coefficients $b$ and $\sigma$ satisfying Assumptions {\bf (H1)}-{\bf (H2)}, take $b(x,i,u)=\bar b(i,u) x$, $\sigma(x,i,u)=\bar\sigma(i,u) x\trans$, for some continuous and bounded $\bar b\colon\I_m\times U\rightarrow\R$, $\bar\sigma\colon\I_m\times U\rightarrow\R^d$ (notice that $\bar\sigma(i,u) x\trans\in\R^{d\times d}$ and $\|\sigma(x,i,u)-\sigma(x',i,u)\|=|\bar\sigma(i,u)||x-x'|$). Then, the dissipative condition \eqref{dissipative} holds if and only if
\begin{equation}\label{Diss_2}
\bar b(i,u) + \frac{1}{2}|\bar\sigma(i,u)|^2 \ \leq \  - \gamma,
\end{equation}
for all $(i,u)\in\I_m\times U$.
\ep
}
\end{Remark}

\subsection{Main results}

Consider the following \emph{ergodic} system of quasi variational inequalities
\begin{equation}
\label{HJB_ergodic}
\min\Big\{\lambda - \inf_{u\in U}\big[\Lc^{i,u}\phi(x,i) + f(x,i,u)\big],\phi(x,i) - \max_{j\neq i}\big[\phi(x,j) - c(x,i,j)\big]\Big\} \ = \ 0.
\end{equation}
We begin providing the definition of viscosity solution to system \eqref{HJB_ergodic}.

\begin{Definition}
\textup{(i)} We say that a pair $(\lambda,\phi)$, with $\lambda\in\R$ and $\phi\colon\R^d\times\I_m\rightarrow\R$ a lower $($resp. upper$)$ semicontinuous function,  is a viscosity supersolution $($resp. subsolution$)$ to the system of variational inequalities \eqref{HJB_ergodic} if
\[
\min\Big\{\lambda - \inf_{u\in U}\big[\Lc^{i,u}\varphi(x) + f(x,i,u)\big],\phi(x,i) - \max_{j\neq i}\big[\phi(x,j) - c(x,i,j)\big]\Big\} \geq (resp.\,\leq) \, 0,
\]
for any $(x,i)\in\R^d\times\I_m$ and any $\varphi\in C^2(\R^d)$ such that
\[
\phi(x,i)-\varphi(x) \ = \ \min_{\R^d}\{\phi(\cdot,i)-\varphi(\cdot)\} \quad \big(\text{resp. }\max_{\R^d}\{\phi(\cdot,i)-\varphi(\cdot)\}\big).
\]
\textup{(ii)} We say that a pair $(\lambda, \phi)$, with $\lambda\in\R$ and $\phi\colon\R^d\times\I_m\rightarrow\R$ a continuous function, is a viscosity solution to the system of variational inequalities \eqref{HJB_ergodic} if it is both a viscosity supersolution and a viscosity subsolution to \eqref{HJB_ergodic}.
\end{Definition}

We can now summarize the main results of the paper in the two following theorems.

\begin{Theorem}
\label{T:ExistHJB_T_beta}
Suppose that Assumption {\bf (H1)} holds.

\noindent\textup{(i)} There exists a unique continuous viscosity solution $V$ to system \eqref{HJBintro} satisfying the growth condition: for any $T>0$, there exist $C_T\geq0$ and $q_T\geq1$ such that
\beqs
|V(t,x,i)| &\leq & C_T\big(1+|x|^{q_T}\big),
\enqs
for all $(t,x,i)\in[0,T]\times\R^d\times\I_m$.

\noindent\textup{(ii)} If in addition Assumption {\bf (H2)} holds, then, for any $\beta>0$, there exists a unique continuous viscosity solution $V^\beta$ to system \eqref{HJB_ellipticintro} satisfying the linear growth condition: there exists $C_\beta\geq0$ such that
\beqs
|V^\beta(x,i)| &\leq & C_\beta\big(1+|x|\big),
\enqs
for all $(x,i)\in\R^d\times\I_m$.
\end{Theorem}

For the next result, we need the following additional assumptions.

\vspace{3mm}

{\bf (H3)} \hspace{3mm} The cost function $c$ is constant with respect to the variable $x\in\R^d$. By an abuse of notation, under Assumption {\bf (H3)}, we denote $c(i,j)=c(x,i,j)$, for all $(x,i,j)\in\R^d\times\I_m\times\I_m$.

\vspace{3mm}

{\bf (H$U$)} \hspace{3mm} The interior set $\mathring U$ of $U$ is connected, and $U$ coincides with the closure of its interior: $U=\text{Cl}(\mathring U)$.

\begin{Remark}
{\rm
(i) Assumption {\bf (H3)} is standard in the literature on switching control pro\-blems. In the present paper, it is used in Propositions \ref{P:V^beta} and \ref{P:VbetaViscSol} (and also in Corollary \ref{C:EstimatesV^beta,n}(i)) to establish the crucial uniform Lipschitz property of $V^\beta$ with respect to $x$.

\noindent(ii) Assumption {\bf (H$U$)} is employed in obtaining a dual representation formula for $V^\beta$ in Proposition \ref{P:VbetaViscSol}.
\ep
}
\end{Remark}

\begin{Theorem}
\label{T:Main}
Let Assumptions {\bf (H1)}, {\bf (H2)}, {\bf (H3)}, and {\bf (H$U$)} hold. Then,  there exists a viscosity solution $(\lambda,\phi)$, with $\phi(\cdot,i)$ Lipschitz, for any $i\in\I_m$, and $\phi(0,i_0)=0$ for some fixed $i_0\in\I_m$, to the ergodic system \eqref{HJB_ergodic}, such that
\[
\beta V^\beta(x,i) \ \overset{\beta\rightarrow0^+}{\longrightarrow} \ \lambda,
\]
for all $(x,i)\in\R^d\times\I_m$, and, for some sequence $(\beta_k)_{k\in\N}$, with $\beta_k\searrow0^+$, we have
\[
V^{\beta_k}(\cdot,i) - V^{\beta_k}(0,i_0) \ \underset{\text{in }C(\R^d)}{\overset{k\rightarrow\infty}{\longrightarrow}} \ \phi(\cdot,i),
\]
for all $i\in\I_m$, where ``\,in $C(\R^d)$\hspace{-.3mm}'' stands for uniform convergence on compact subsets of $\R^d$. Moreover, for any viscosity solution $(\lambda,\phi)$ to \eqref{HJB_ergodic}, with $\phi$ satisfying
\begin{equation}
\label{GrowthCond_phi2_Main}
|\phi(x,i)| \ \leq \ M_\phi(1 + |x|^2), \qquad \forall\,(x,i)\in\R^d\times\I_m,
\end{equation}
for some constants $M_\phi\geq0$, we have
\[
\frac{V(T,x,i)}{T} \ \overset{T\rightarrow\infty}{\longrightarrow} \ \lambda, \qquad \forall\,(x,i)\in\R^d\times\I_m.
\]
In particular, $\lambda$ is uniquely determined for all viscosity solutions $(\lambda,\phi)$ to \eqref{HJB_ergodic}, with $\phi$ satisfying a quadratic growth condition as in \eqref{GrowthCond_phi2_Main}.
\end{Theorem}

\begin{Remark}
{\rm
{\bf (i)} A question which naturally arises from Theorem \ref{T:Main} is the uniqueness of $\phi$. This problem has been tackled in \cite{CFP14}, Theorem 5.2, using probabilistic techniques as in \cite{ichihara12} and assuming some smoothness of $\phi$ together with the existence of an optimal feedback control (under which the state process admits a unique invariant measure). This kind of proof seems, however, designed to deal with ergodic equation associated to stochastic control problems rather than to stochastic differential games. Let us explain more in detail this latter issue, recalling the main steps of the proof of Theorem 5.2 in \cite{CFP14} (translated into the present framework) and emphasizing where it breaks down for stochastic differential games. The goal in \cite{CFP14} consists in proving, for every $i\in\I_m$, that there exists a real constant $C_i$ such that
\begin{equation}
\label{W-->C_i}
W(T,x,i) \ := \ V(T,x,i) - \big(\lambda T + \phi(x,i)\big) \ \overset{T\rightarrow\infty}{\longrightarrow} \ C_i, \qquad \forall\,(x,i)\in\R^d\times\I_m.
\end{equation}
We see that if \eqref{W-->C_i} holds for every viscosity solution $(\lambda,\phi)$, with $\phi$ Lipschitz, to the ergodic system \eqref{HJB_ergodic}, then, for every $i\in\I_m$, $\phi(\cdot,i)$ is uniquely determined up to a constant which depends only on $i\in\I_m$. To prove \eqref{W-->C_i}, we proceed as follows (we just sketch the main steps and do not pause on the technicalities, however we use some results from Section \ref{S:Parabolic} below). For any $T,S\geq0$, from the identity $V(T+S,x,i)=V^{T+S}(0,x,i,)$ stated in \eqref{V=V^T}, we deduce
\begin{align*}
V(T+S,x,i) \ &= \sup_{\alpha\in\Ac_{0,0}}\inf_{\upsilon\in\Uc_{0,0}} \E\bigg[\int_0^T f(X_s^{0,x,i;\alpha,\upsilon},I_s^{0,x,i;\alpha,\upsilon},\upsilon_s)ds + V(S,X_T^{0,x,i;\alpha,\upsilon},I_T^{0,x,i;\alpha,\upsilon}) \\
&\quad - \sum_{n\in\N} c(X_{\tau_n}^{0,x,i;\alpha,\upsilon},I_{\tau_n^-}^{0,x,i;\alpha,\upsilon},I_{\tau_n}^{0,x,i;\alpha,\upsilon}) 1_{\{\tau_n<T\}} \bigg].
\end{align*}
On the other hand, from \eqref{phi=J} with $t=0$, we have
\begin{align*}
\phi(x,i) \ &= \sup_{\alpha\in\Ac_{0,0}}\inf_{\upsilon\in\Uc_{0,0}} \E\bigg[\int_0^T \big(f(X_s^{0,x,i;\alpha,\upsilon},I_s^{0,x,i;\alpha,\upsilon},\upsilon_s) - \lambda\big)ds + \phi(X_T^{0,x,i;\alpha,\upsilon},I_T^{0,x,i;\alpha,\upsilon}) \\
&\quad - \sum_{n\in\N} c(X_{\tau_n}^{0,x,i;\alpha,\upsilon},I_{\tau_n^-}^{0,x,i;\alpha,\upsilon},I_{\tau_n}^{0,x,i;\alpha,\upsilon}) 1_{\{\tau_n<T\}} \bigg].
\end{align*}
Now, suppose that there exists an optimal control $\alpha^*=(\tau_n^*,\iota_n^*)_{n\in\N}\in\Ac_{0,0}$, independent of $T$, such that
\begin{align*}
\phi(x,i) \ &= \ \inf_{\upsilon\in\Uc_{0,0}} \E\bigg[\int_0^T \big(f(X_s^{0,x,i;\alpha^*,\upsilon},I_s^{0,x,i;\alpha^*,\upsilon},\upsilon_s) - \lambda\big)ds + \phi(X_T^{0,x,i;\alpha^*,\upsilon},I_T^{0,x,i;\alpha^*,\upsilon}) \\
&\quad \ - \sum_{n\in\N} c(X_{\tau_n^*}^{0,x,i;\alpha^*,\upsilon},I_{(\tau_n^*)^-}^{0,x,i;\alpha^*,\upsilon},I_{\tau_n^*}^{0,x,i;\alpha^*,\upsilon}) 1_{\{\tau_n^*<T\}} \bigg].
\end{align*}
As a consequence, we obtain the inequality
\begin{equation}
\label{W_T+SgeqW_T}
W(T+S,x,i) \ \geq \ \inf_{\upsilon\in\Uc_{0,0}} \E\big[W(S,X_T^{0,x,i;\alpha^*,\upsilon},I_T^{0,x,i;\alpha^*,\upsilon})\big].
\end{equation}
As in \cite{CFP14} and \cite{ichihara12}, let us introduce, for any $i\in\I_m$, the set of $\omega$-limits of $\{W(T,\cdot,i)\}_{T\geq0}$ in $C(\R^d)$:
\[
\Gamma_i \ := \ \big\{w_\infty\in C(\R^d)\colon W(T_j,\cdot,i)\rightarrow w_\infty\text{ in }C(\R^d)\text{ for some }T_j\nearrow\infty\big\}.
\]
It can be shown that $\Gamma_i\neq\emptyset$. Therefore, \eqref{W-->C_i} follows if we prove that $\Gamma_i=\{w_{\infty,i}\}$ is a singleton, with $w_{\infty,i}$ constant. In \cite{CFP14}, Theorem 5.2, the idea is to pass to the limit in \eqref{W_T+SgeqW_T}, using that under $\alpha^*$ the state process admits a unique invariant measure. Since, however, in \eqref{W_T+SgeqW_T} there is also the ``inf'' operator (due to the \emph{game} feature of the robust switching control problem), the same argument does not allow to conclude (indeed, e.g., the ``invariant measure'' would depend on $\upsilon\in\Uc_{0,0}$).

\noindent {\bf (ii)} Notice that $\lambda$ appearing in Theorem \ref{T:Main} does not depend on $i\in\I_m$.

\noindent {\bf (iii)} Assumptions {\bf (H3)} and {\bf (H$U$)} in Theorem \ref{T:Main} are not needed to prove the convergence results of $V$, but only those of $V^\beta$.

\noindent {\bf (iv)} Theorem \ref{T:Main} can be interpreted as a \emph{Tauberian theorem}. Indeed, the convergence results for $V^\beta$ allows to prove the existence of a viscosity solution $(\lambda,\phi)$ to \eqref{HJB_ergodic}, from which the asymptotic results for $V$ follow (see also Remark 5.3(ii) in \cite{CFP14}).
\ep
}
\end{Remark}

\begin{Remark} \label{remergodic}
{\rm
Suppose that Assumptions {\bf (H1)}, {\bf (H2)}, {\bf (H3)}, and {\bf (H$U$)} hold. Suppose also that there exists a viscosity solution $(\lambda,\phi)$ to \eqref{HJB_ergodic} with $\phi$ satisfying \eqref{GrowthCond_phi2_Main}. Finally, similarly to \eqref{Phi_alpha*}, suppose that for every $\eps>0$ there exists an $\eps$-optimal control (we omit the dependence on $\eps$) $\alpha^*=(\tau_n^*,\iota_n^*)_{n\in\N}\in\Ac_{0,0}$ (we refer to Section \ref{SubS:Robust} for all unexplained notations), independent of $T$, such that
\begin{align}
\phi(x,i) \ &\leq \ \inf_{\upsilon\in\Uc_{0,0}} \E\bigg[\int_0^T \big(f(X_s^{0,x,i;\alpha^*,\upsilon},I_s^{0,x,i;\alpha^*,\upsilon},\upsilon_s) - \lambda\big)ds + \phi(X_T^{0,x,i;\alpha^*,\upsilon},I_T^{0,x,i;\alpha^*,\upsilon}) \notag \\
&\quad \ - \sum_{n\in\N} c(X_{\tau_n^*}^{0,x,i;\alpha^*,\upsilon},I_{(\tau_n^*)^-}^{0,x,i;\alpha^*,\upsilon},I_{\tau_n^*}^{0,x,i;\alpha^*,\upsilon}) 1_{\{\tau_n^*<T\}} \bigg] + \eps. \label{Phi_alpha*}
\end{align}
Then $\lambda$ can be interpreted as value of a \emph{robust ergodic control problem}:
\begin{equation}\label{lambda=J}
\lambda \ = \ \sup_{\alpha\in\Ac_{0,0}} J(x,i;\alpha), \qquad \forall\,(x,i)\in\R^d\times\I_m,
\end{equation}
with
\begin{align*}
J(x,i;\alpha) \ &:= \ \limsup_{T\rightarrow\infty}\frac{1}{T}\inf_{\upsilon\in\Uc_{0,0}}\E\bigg[\int_0^T f(X_t^{0,x,i;\alpha,\upsilon},I_t^{0,x,i;\alpha,\upsilon},\upsilon_t)dt \\
&\quad \ \ - \sum_{n\in\N} c(X_{\tau_n}^{0,x,i;\alpha,\upsilon},I_{\tau_n^-}^{0,x,i;\alpha,\upsilon},I_{\tau_n}^{0,x,i;\alpha,\upsilon}) 1_{\{\tau_n<T\}} \bigg],
\end{align*}
where $\tau^n$ stands for $\tau^n(X_\cdot^{0,x,i;\alpha,\upsilon},I_{\cdot^-}^{0,x,i;\alpha,\upsilon})$, and the state processes $X^{0,x,i;\alpha,\upsilon},I^{0,x,i;\alpha,\upsilon}$ satisfy system \eqref{SDE} below, with $t=0$, $(x,i)\in\R^d\times\I_m$, $\alpha=(\tau_n,\iota_n)_{n\in\N}\in\Ac_{0,0}$.

Let us prove \eqref{lambda=J}. Firstly, observe that the viscosity solution $V$ to system \eqref{HJBintro}, whose existence is stated in Theorem \ref{T:ExistHJB_T_beta}, admits the stochastic control representation (see identity \eqref{V=V^T} below):
\[
V(T,x,i) \ = \ V^T(0,x,i) \ = \ \sup_{\alpha\in\Ac_{0,0}}\inf_{\upsilon\in\Uc_{0,0}} J^T(0,x,i;\alpha,\upsilon),
\]
where
\begin{align*}
J^T(0,x,i;\alpha,\upsilon) \ &= \ \E\bigg[\int_0^T f(X_t^{0,x,i;\alpha,\upsilon},I_t^{0,x,i;\alpha,\upsilon},\upsilon_t)dt + g(X_T^{0,x,i;\alpha,\upsilon},I_T^{0,x,i;\alpha,\upsilon}) \\
&\quad \ - \sum_{n\in\N} c(X_{\tau_n}^{0,x,i;\alpha,\upsilon},I_{\tau_n^-}^{0,x,i;\alpha,\upsilon},I_{\tau_n}^{0,x,i;\alpha,\upsilon}) 1_{\{\tau_n<T\}} \bigg].
\end{align*}
Then, by Theorem \ref{T:Main}, for any $(x,i)\in\R^d\times\I_m$,
\begin{align}
\lambda \ &= \ \lim_{T\rightarrow\infty}\frac{V(T,x,i)}{T} \ = \ \lim_{T\rightarrow\infty}\sup_{\alpha\in\Ac_{0,0}}\inf_{\upsilon\in\Uc_{0,0}}\frac{1}{T}\E\bigg[ \int_0^T f(X_t^{0,x,i;\alpha,\upsilon},I_t^{0,x,i;\alpha,\upsilon},\upsilon_t) dt \notag \\
&\quad \ + g(X_T^{0,x,i;\alpha,\upsilon},I_T^{0,x,i;\alpha,\upsilon}) - \sum_{n\in\N} c(X_{\tau_n}^{0,x,i;\alpha,\upsilon},I_{\tau_n^-}^{0,x,i;\alpha,\upsilon},I_{\tau_n}^{0,x,i;\alpha,\upsilon}) 1_{\{\tau_n<T\}} \bigg] \notag \\
&= \ \lim_{T\rightarrow\infty}\sup_{\alpha\in\Ac_{0,0}}\inf_{\upsilon\in\Uc_{0,0}}\frac{1}{T}\E\bigg[ \int_0^T f(X_t^{0,x,i;\alpha,\upsilon},I_t^{0,x,i;\alpha,\upsilon},\upsilon_t) dt \notag \\
&\quad \ - \sum_{n\in\N} c(X_{\tau_n}^{0,x,i;\alpha,\upsilon},I_{\tau_n^-}^{0,x,i;\alpha,\upsilon},I_{\tau_n}^{0,x,i;\alpha,\upsilon}) 1_{\{\tau_n<T\}} \bigg], \label{lambda=lambda^*}
\end{align}
where the last equality follows from the fact that
\[
\lim_{T\rightarrow\infty}\sup_{\alpha\in\Ac_{0,0}}\inf_{\upsilon\in\Uc_{0,0}}\frac{1}{T}\E\big[g(X_T^{0,x,i;\alpha,\upsilon},I_T^{0,x,i;\alpha,\upsilon})\big] \ = \ 0,
\]
which is a consequence of the quadratic growth condition of $g$ and estimate \eqref{EstimateX_diss}.
From \eqref{lambda=lambda^*} we see that $\lambda\geq\sup_{\alpha\in\Ac_{0,0}} J(x,i;\alpha)$. To prove the reverse inequality, fix $(x,i)\in\R^d\times I_m$ and $\eps>0$, then by \eqref{Phi_alpha*} we obtain
\begin{align*}
\lambda \ &\leq \ \frac{1}{T} \inf_{\upsilon\in\Uc_{0,0}} \E\bigg[\int_0^T f(X_t^{0,x,i;\alpha^*,\upsilon},I_t^{0,x,i;\alpha^*,\upsilon},\upsilon_t)dt + \phi(X_T^{0,x,i;\alpha^*,\upsilon},I_T^{0,x,i;\alpha^*,\upsilon}) - \phi(x,i) \\
&\quad \ - \sum_{n\in\N} c(X_{\tau_n^*}^{0,x,i;\alpha^*,\upsilon},I_{(\tau_n^*)^-}^{0,x,i;\alpha^*,\upsilon},I_{\tau_n^*}^{0,x,i;\alpha^*,\upsilon}) 1_{\{\tau_n^*<T\}} \bigg] + \eps.
\end{align*}
From the Lipschitz property of $\phi$ and estimate \eqref{EstimateX_diss}, we have
\[
\lim_{T\rightarrow\infty}\frac{1}{T} \inf_{\upsilon\in\Uc_{0,0}} \E\big[\phi(X_T^{0,x,i;\alpha^*,\upsilon},I_T^{0,x,i;\alpha^*,\upsilon}) - \phi(x,i)\big] \ = \ 0,
\]
therefore
\begin{align*}
\lambda \ &\leq \ \lim_{T\rightarrow\infty} \frac{1}{T} \inf_{\upsilon\in\Uc_{0,0}} \E\bigg[\int_0^T f(X_t^{0,x,i;\alpha^*,\upsilon},I_t^{0,x,i;\alpha^*,\upsilon},\upsilon_t)dt \\
&\quad \ - \sum_{n\in\N} c(X_{\tau_n^*}^{0,x,i;\alpha^*,\upsilon},I_{(\tau_n^*)^-}^{0,x,i;\alpha^*,\upsilon},I_{\tau_n^*}^{0,x,i;\alpha^*,\upsilon}) 1_{\{\tau_n^*<T\}} \bigg] + \eps \ \leq \ \sup_{\alpha\in\Ac_{0,0}}J(x,i;\alpha) + \eps.
\end{align*}
From the arbitrariness of $\eps$, we find $\lambda\leq \sup_{\alpha\in\Ac_{0,0}}J(x,i;\alpha)$, which, together with \eqref{lambda=lambda^*}, yields \eqref{lambda=J}.
\ep
}
\end{Remark}

\begin{Remark}
{\rm
\emph{Explicit solution to \eqref{HJB_ergodic} in the two-regime case.} We show the validity of Theorem \ref{T:Main} in a specific example, where we are able to find an explicit solution to the ergodic system of quasi variational inequalities \eqref{HJB_ergodic}. More precisely, we consider the framework of Section 5.3 in \cite{pham09}, where an infinite horizon two-regime switching control problem is studied and an explicit solution is determined. We slightly generalize the results in \cite{pham09} in order to take into account the robust feature. This allows us to find explicitly $V^\beta$, $\beta>0$. Then, as expected from Theorem \ref{T:Main}, letting $\beta\rightarrow0^+$, we are able to construct a solution to \eqref{HJB_ergodic}.

Take $d=1$ and $m=2$. Let $b$ and $\sigma$ be as in Remark \ref{R:coefficients}.(ii), so that the state process $X$ can assume only positive values whenever the initial condition is positive. Let also $c(x,i,j)=c(i,j)$ be independent of $x$, and $f(x,i,u)=x^p$, $x>0$, for some $p\in(0,1)$. 

Define the constants $\kappa_1$ and $\kappa_2$ as
\[
\kappa_i \ = \ - \inf_{u\in U}\big[\bar b(i,u)p + \tfrac{1}{2}\bar\sigma(i,u)^2 p (p-1)\big], \qquad i=1,2.
\]
Recalling from condition \eqref{Diss_2} that $\bar b$ is a strictly negative function, we see that $\kappa_i>0$, $i=1,2$. Proceeding as in Theorem 5.3.4 of \cite{pham09}, we now find the explicit expression of $V^\beta$, after distinguishing between the two cases: $\kappa_1=\kappa_2$ and $\kappa_1\neq\kappa_2$ (in this latter case we need to impose an additional assumption on $\bar b$ and $\bar\sigma$ in order to follow the same steps as in the proof of Theorem 5.3.4).

\vspace{1mm}

\noindent\emph{Case 1:} $\kappa_1=\kappa_2$. Consider the function
\[
V^\beta(x,i) \ = \ K_i^\beta x^p, \qquad x>0,\;i=1,2,
\]
where
\[
K_i^\beta \ = \ \frac{1}{\beta + \kappa_i}, \qquad i=1,2.
\]
By direct calculation, we can prove that $V^\beta$ is a viscosity solution to equation \eqref{HJB_ellipticintro} on $(0,\infty)\times\{1,2\}$. Notice that this latter result is still true when $\beta=0$. Then, it is easy to see that the pair $(\lambda,\phi)$ given by $\lambda=\lim_{\beta\rightarrow0^+}\beta V^\beta(x,i)=0$ and $\phi(\cdot,i)=\lim_{\beta\rightarrow0^+}V^\beta(\cdot,i)=K_i^0 {\cdot\,}^p$, with $K_i^0=1/\kappa_i$, is a viscosity solution to the ergodic system \eqref{HJB_ergodic}. From the explicit expression of the pair $(\lambda,\phi)$, proceeding as in Theorem 5.3.4 in \cite{pham09}, we can also determine an optimal switching control for the robust ergodic control problem \eqref{lambda=J}. Indeed, it is easy to see that, in this case, it is never optimal to switch.

\vspace{1mm}

\noindent\emph{Case 2:} $\kappa_1\neq\kappa_2$. In order to reason as in the proof of Theorem 5.3.4 in \cite{pham09}, we impose the following additional assumption: $\bar b(i,u)=\hat b(i)h(u)$ and $\bar\sigma(i,u)=\hat\sigma(i)\sqrt{h(u)}$, with $\hat b\colon\{1,2\}\rightarrow\R$, $\hat\sigma\colon\{1,2\}\rightarrow\R$, $h\colon U\rightarrow(0,\infty)$ continuous and bounded, satisfying (see condition \eqref{Diss_2})
\[
\hat b(i)h(u) + \frac{1}{2}\hat\sigma(i)^2 h(u) \ \leq \ - \gamma,
\]
for all $(i,u)\in\{1,2\}\times U$. In particular, $\hat b$ is a strictly negative function.

Without loss of generality, we suppose that $\kappa_1>\kappa_2$, the other case can be treated in an analogous way. Let
\[
m_\beta \ = \ \frac{1}{2} - \frac{\hat b(1)}{\hat\sigma(1)^2} + \sqrt{\bigg(\frac{1}{2} - \frac{\hat b(1)}{\hat\sigma(1)^2}\bigg)^2 + \frac{2\beta}{\hat\sigma(1)^2\inf_{u\in U}h(u)}} \ > \ 1.
\]
Then, proceeding as in Theorem 5.3.4 of \cite{pham09}, we see that the unique continuous viscosity solution to equation \eqref{HJB_ellipticintro} on $(0,\infty)\times\{1,2\}$ is given by the function $V^\beta$ defined as follows:
\begin{align*}
V^\beta(x,1) \ &= \
\begin{cases}
A_\beta x^{m_\beta} + K_1^\beta x^p, \qquad & x\in(0,x_\beta), \\
K_2^\beta x^p - c(1,2), & x\in[x_\beta,\infty),
\end{cases} \\
V^\beta(x,2) \ &= \ K_2^\beta x^p, \hspace{2.3cm} x\in(0,\infty),
\end{align*}
with
\begin{align*}
x_\beta \ &= \ \bigg(\frac{m_\beta}{m_\beta-p}\frac{c(1,2)}{K_2^\beta-K_1^\beta}\bigg)^{\frac{1}{p}}, \\
A_\beta \ &= \ \big(K_2^\beta - K_1^\beta\big) \frac{p}{m_\beta} x_\beta^{p-m_\beta}.
\end{align*}
Notice that $x_\beta$ and $A_\beta$ are determined by the continuity and smooth-pasting conditions of $V^\beta(\cdot,1)$ at $x_\beta$:
\begin{align*}
A_\beta x_\beta^{m_\beta} + K_1^\beta x_\beta^p \ &= \ K_2^\beta x_\beta^p - c(1,2), \\
A_\beta m_\beta x_\beta^{m_\beta-1} + K_1^\beta p x_\beta^{p-1} \ &= \ K_2^\beta p x_\beta^{p-1}.
\end{align*}
Then, it is easy to see that a viscosity solution to the ergodic system \eqref{HJB_ergodic} is given by the pair $(\lambda,\phi)$, with $\lambda=\lim_{\beta\rightarrow0^+}\beta V^\beta(x,i)=0$ and $\phi(\cdot,i)=\lim_{\beta\rightarrow0^+} V^\beta(\cdot,i)$ defined as
\begin{align*}
\phi(x,1) \ &= \ \begin{cases}
A_0 x^{m_0} + K_1^0 x^p, \qquad & x\in(0,x_0), \\
K_2^0 x^p - c(1,2), & x\in[x_0,\infty),
\end{cases} \\
\phi(\cdot,2) \ &= \ K_2^0 x^p, \hspace{2.3cm} x\in(0,\infty),
\end{align*}
where $K_1^0$, $K_2^0$, $m_0$, $x_0$, $A_0$ correspond to $K_1^\beta$, $K_2^\beta$, $m_\beta$, $x_\beta$, $A_\beta$ with $\beta=0$. Exploiting the knowledge of $(\lambda,\phi)$, we can also find an optimal switching control for the robust ergodic control problem \eqref{lambda=J}. Indeed, proceeding as in Theorem 5.3.4 in \cite{pham09}, it is easy to see that: when we are in regime $1$, it is optimal to switch to regime $2$ whenever the state process $X$ exceeds the threshold $x_0$; while it is never optimal to switch when we are in regime $2$.
\ep}
\end{Remark}

The rest of the paper is devoted not only to the proof of Theorems \ref{T:ExistHJB_T_beta} and \ref{T:Main}, but also to investigate more in detail the properties of the systems \eqref{HJBintro}, \eqref{HJB_ellipticintro}, and \eqref{HJB_ergodic}, in particular from a stochastic control point of view, exploring their relation with robust switching control problems. To sum up, the logical flow of the paper is the following: Section \ref{S:Parabolic} is devoted to the analysis of the parabolic system. The essence of this section are three propositions: Proposition 
\ref{P:ExistV} states the connection between the parabolic system of variational inequalities \reff{HJBintro} and robust switching control problem. 
Proposition \ref{P:repphi} gives robust control representation bounds for a solution to the ergodic system, which together with Proposition \ref{P:ExistV}  leads to a quite simple argument for proving in Proposition \ref{P:V-->lambda} the long time convergence of the parabolic system.  In Section \ref{S:Elliptic}, we analyze the elliptic system. The essence of this section can be summarized as follows.
Proposition  \ref{T:Feynman-Kac} provides a Feynman-Kac formula in terms of BSDE for the solution $V^{\beta,n}$ to the penalized system of elliptic variational inequalities, which is then used for stating in Corollary \ref{C:DualFormula}, a dual representation for $V^{\beta,n}$, and then by passing to the limit in $n$, for getting in Proposition \ref{P:V^beta}  a dual probabilistic game representation for $V^\beta$.  From these dual representations, we are able to derive key uniform estimates for $V^{\beta,n}$ in Corollary \ref{C:EstimatesV^beta,n}, and then  for $V^\beta$ in Propositions \ref{P:V^beta} and 
\ref{P:VbetaViscSol}. Finally, we can obtain the existence of a solution to the ergodic system and the convergence result of $V$ and $V^\beta$ in Proposition 
\ref{P:convergence}.

\section{Long time asymptotics of the parabolic system}
\label{S:Parabolic}

\setcounter{equation}{0} 
\setcounter{Theorem}{0} \setcounter{Proposition}{0}
\setcounter{Corollary}{0} \setcounter{Lemma}{0}
\setcounter{Definition}{0} \setcounter{Remark}{0}

We firstly investigate the long time asymptotics of the solution $V$ to the parabolic system of variational inequalities \eqref{HJBintro}. To this end, we shall rely on probabilistic arguments, which are based on a characterization of $V$ in terms of a finite horizon robust switching control problem introduced in \cite{BCP14}, that we now describe.

\subsection{Robust feedback switching control problem}
\label{SubS:Robust}

We present the robust switching control problem focusing only on the main issues, in order to alleviate the presentation and reduce as much as possible the technicalities, for which we refer to \cite{BCP14}. 

\vspace{1mm}

Consider a complete probability space $(\Omega,\Fc,\P)$ and a $d$-dimensional Brownian motion $W=(W_t)_{t\geq0}$ defined on it. Denote by $\F=(\Fc_t)_{t\geq0}$ the completion of the natural filtration generated by $W$. Fix a finite time horizon $T\in(0,\infty)$. In the robust switching control problem we are going to present, the switcher plays against an adverse player, which can be interpreted as nature and renders the optimization problem robust. We begin recalling the type of controls adopted by the switcher and by nature, following Definitions 2.1 and 2.2 in \cite{BCP14}, to which we refer for more details:
\begin{itemize}
\item $\Ac_{t,t}$, $t\in[0,T]$, denotes the family of all \emph{feedback switching controls} starting at time $t$ for the switcher. A generic element of $\Ac_{t,t}$ is given by a double sequence $\alpha=(\tau_n,\iota_n)_{n\in\N}$, where $(\tau_n)_{n\in\N}$ is a nondecreasing sequence of stopping times valued in $[t,T]$ and $\iota_n\in\I_m$ represents the switching action, i.e., the regime from time $\tau_n$ up to time $\tau_{n+1}$. The random variable $\iota_n$ only depends on the information known up to time $\tau_n$. Each $\alpha$ has a feedback form, in the sense that it is chosen  by the switcher based only on the past and present information coming from the state and regime processes.
\item $\Uc_{t,t}$, $t\in[0,T]$, denotes the family of all \emph{open-loop controls} starting at time $t$ for nature. A generic element of $\Uc_{t,t}$ is an adapted process $\upsilon\colon[t,T]\times\Omega\rightarrow U$.
\end{itemize}
As explained in \cite{BCP14}, the reason behind the feedback form of a switching control comes from the observation that the switcher in general knows only the evolution of the state process $X$ and regime $I$. On the other hand, a control $\upsilon\in\Uc_{t,t}$ is not necessarily of feedback form, since nature at time $t\in[0,T]$ is wise to all information up to time $t$.

\vspace{1mm}

Let us now introduce the controlled dynamics of the state and regime processes. For any $(t,x,i)\in[0,T]\times\R^d\times\I_m$, $\alpha=(\tau_n,\iota_n)_{n\in\N}\in\Ac_{t,t}$, $\upsilon\in\Uc_{t,t}$, the state process $X$ and regime $I$ evolve on $[t,T]$ according to the following controlled SDEs:
\begin{equation}
\label{SDE}
\begin{cases}
X_s \ \!\!= \ x + \int_t^s b(X_r,I_r,\upsilon_r)dr + \int_t^s\sigma(X_r,I_r,\upsilon_r)dW_r, &t\leq s\leq T, \\
I_s \ = \ i1_{\{t\leq s<\tau_0(X_\cdot,I_{\cdot^-})\}} + \sum_{n\in\N} \iota_n(X_\cdot,I_{\cdot^-}) 1_{\{\tau_n(X_\cdot,I_{\cdot^-})\leq s<\tau_{n+1}(X_\cdot,I_{\cdot^-})\}}, &t\leq s<T, \\
I_T \ \!\!= \ I_{T^-},
\end{cases}
\end{equation}
with $I_{t^-}:=I_t$. Notice that $\tau_n$ and $\iota_n$ have a feedback form, indeed they depend only on $X$ and $I$. We recall the following wellposedness result from \cite{BCP14}.

\begin{Lemma}
\label{L:X}
Suppose that Assumption {\bf (H1)} holds. For any $T>0$, $(t,x,i)\in[0,T]\times\R^d\times\I_m$, $\alpha\in\Ac_{t,t}$, $\upsilon\in\Uc_{t,t}$, there exists a unique (up to indistinguishability) $\F$-adapted process $(X^{t,x,i;\alpha,\upsilon},I^{t,x,i;\alpha,\upsilon})=(X_s^{t,x,i;\alpha,\upsilon},I_s^{t,x,i;\alpha,\upsilon})_{t\leq s\leq T}$ to equation \eqref{SDE}. Moreover, for any $p\geq2$ there exists a positive constant $C_{p,T}$, depending only on $p,T,L_1$ (independent of $t,x,i,\alpha,\upsilon$), such that
\beq
\label{EstimateX}
\E\Big[\sup_{t\leq s\leq T}|X_s^{t,x,i;\alpha,\upsilon}|^p\Big] &\leq & C_{p,T}(1+|x|^p).
\enq
\end{Lemma}
\textbf{Proof.}
See Proposition 2.1 in \cite{BCP14}.
\ep

\vspace{3mm}

We also have the following result as a consequence of the dissipativity condition.

\begin{Lemma}
\label{L:X_diss}
Suppose that Assumptions {\bf (H1)} and {\bf (H2)} hold. There exists a positive constant $\bar C$, depending only on $M_{b,\sigma}:=\sup_{(i,u)\in\I_m\times U}(|b(0,i,u)|+\|\sigma(0,i,u)\|)$, $L_1$, and $\gamma$, such that
\beq
\label{EstimateX_diss}
\sup_{s\in[t,T],\,\alpha\in\Ac_{t,t},\,\upsilon\in\Uc_{t,t}}\E\big[|X_s^{t,x,i;\alpha,\upsilon}|^2\big] &\leq & \bar C(1+|x|^2),
\enq
for any $(t,x,i)\in[0,T]\times\R^d\times\I_m$ and $T>0$.
\end{Lemma}
\textbf{Proof.}
Fix $T>0$, $(t,x,i)\in[0,T]\times\R^d\times\I_m$, $\alpha\in\Ac_{t,t}$, $\upsilon\in\Uc_{t,t}$. The proof can be done along the lines of Lemma 2.1(i) in \cite{CFP14}. We simply recall the main steps. We take $s\in[t,T]$ and apply It\^o's formula to $e^{\gamma(r-t)}|X_r^{t,x,i;\alpha,\upsilon}|^2$ between $r=t$ and $r=s$. Then, we rearrange the terms in order to exploit the dissipativity condition {\bf (H2)}. Afterwards, using the uniform linear growth condition of $b$ and $\sigma$ with respect to $x$, we find that there exists a constant $\bar C$, depending only on $M_{b,\sigma}$, $L_1$, $\gamma$, such that
\begin{equation}
\label{E:ItoProof2}
|X_s^{t,x,i;\alpha,\upsilon}|^2 \ \leq \ \bar C\Big(1 + |x|^2 + \int_t^s e^{\gamma(r-s)}(X_r^{t,x,i;\alpha,\upsilon})\trans\sigma(X_r^{t,x,i;\alpha,\upsilon},I_r^{t,x,i;\alpha,\upsilon},\upsilon_r)dW_r\Big).
\end{equation}
From estimate \eqref{EstimateX} and the linear growth of $\sigma$, we see that the stochastic integral in \eqref{E:ItoProof2} is a martingale. Therefore, taking the expectation in \eqref{E:ItoProof2}, the claim follows.
\ep

\vspace{3mm}

We can now introduce the value function of the robust feedback switching control problem, which is given by:
\beq
\label{V^T}
V^T(t,x,i) &:=& \sup_{\alpha\in\Ac_{t,t}}\inf_{\upsilon\in\Uc_{t,t}} J^T(t,x,i;\alpha,\upsilon), \qquad \forall\,(t,x,i)\in[0,T]\times\R^d\times\I_m,
\enq
with
\beqs
J^T(t,x,i;\alpha,\upsilon) &:=& \E\bigg[\int_t^T f(X_s^{t,x,i;\alpha,\upsilon},I_s^{t,x,i;\alpha,\upsilon},\upsilon_s)ds + g(X_T^{t,x,i;\alpha,\upsilon},I_T^{t,x,i;\alpha,\upsilon}) \\
& & \quad - \; \sum_{n\in\N} c(X_{\tau_n}^{t,x,i;\alpha,\upsilon},I_{\tau_n^-}^{t,x,i;\alpha,\upsilon},I_{\tau_n}^{t,x,i;\alpha,\upsilon}) 1_{\{\tau_n<T\}} \bigg],
\enqs
where $\tau^n$ stands for $\tau^n(X_\cdot^{s,x,i;\alpha,\upsilon},I_{\cdot^-}^{s,x,i;\alpha,\upsilon})$. Notice that the presence of the $\inf_{\upsilon\in\Uc_{t,t}}$ in \eqref{V^T} means that we are looking at the worst case scenario for the switcher and makes the switching control problem robust.

\vspace{1mm}

The dynamic programming equation associated to the robust switching control problem is given by the following system of backward parabolic variational inequalities:
\beq
\label{HJB_T}
\begin{cases}
\min\Big\{- \dfrac{\partial V^T}{\partial t}(t,x,i) - \inf_{u\in U}\big[\Lc^{i,u}V^T(t,x,i) + f(x,i,u)\big], \\
V^T(t,x,i) - \max_{j\neq i}\big[V^T(t,x,j) - c(x,i,j)\big]\Big\} \ = \ 0, &\!\!\!\!\!\!\!\!(t,x,i)\in[0,T)\times\R^d\times\I_m, \\
V^T(T,x,i) \ = \ g(x,i), &\!\!\!\!\!\!\!\!(x,i)\in\R^d\times\I_m.
\end{cases}
\enq
In Corollary 4.1 of \cite{BCP14} it is proved, by means of the stochastic Perron method, that $V^T$ satisfies the dynamic programming principle and it is the unique continuous viscosity solution to system \eqref{HJB_T} (see Definition 2.3 in \cite{BCP14} for the definition of viscosity solution to \eqref{HJB_T}), satisfying a polynomial growth condition
\begin{equation}
\label{GrowthCond_V^T}
\sup_{(t,x,i)\in[0,T]\times\R^d\times\I_m}\frac{|V^T(t,x,i)|}{1+|x|^q} \ < \ \infty,
\end{equation}
for some $q\geq1$.

\vspace{3mm}

We can now present the relation between system \eqref{HJBintro} and the robust switching control problem, which also gives a wellposedness result for viscosity solutions to \eqref{HJBintro} (we do not recall here the definition of viscosity solution to \eqref{HJBintro}, since it is standard and similar to the definition of viscosity solution to \eqref{HJB_T}, for which we refer to Definition 2.3 in \cite{BCP14}).

\begin{Proposition}
\label{P:ExistV}
Suppose that Assumption {\bf (H1)} holds. Then, there exists a unique continuous viscosity solution $V$ to system \eqref{HJBintro} satisfying the growth condition: for any $T>0$, there exist $C_T\geq0$ and $q_T\geq1$ such that
\beq
\label{GrowthCond_T}
|V(t,x,i)| &\leq & C_T\big(1+|x|^{q_T}\big),
\enq
for all $(t,x,i)\in[0,T]\times\R^d\times\I_m$. The function $V$ is given by
\begin{equation}
\label{V=V^T}
V(t,x,i) \ := \ V^T(T-t,x,i), \qquad \forall\,(t,x,i)\in[0,T]\times\R^d\times\I_m,
\end{equation}
for any $T>0$, where $V^T$ is defined by \eqref{V^T}.
\end{Proposition}
\begin{Remark}
{\rm
Notice that point (i) of Theorem \ref{T:ExistHJB_T_beta} follows from Proposition \ref{P:ExistV}.
\ep
}
\end{Remark}
\textbf{Proof.}
\textbf{Step I.} \emph{Existence.} We begin noting that
\begin{equation}
\label{Identification_T0}
V^T(s,x,i) \ = \ V^{T'}(s+T'-T,x,i), \qquad \forall\,(s,x,i)\in[0,T]\times\R^d\times\I_m,
\end{equation}
for any $0\leq T\leq T'<\infty$. Indeed, $V^{T'}(\cdot+T'-T,\cdot,\cdot)$ is a viscosity solution to \eqref{HJB_T} on $[0,T]\times\R^d\times\I_m$, so that identification \eqref{Identification_T0} follows from comparison Theorem 4.1 in \cite{BCP14}. Setting $t:=T-s$ in \eqref{Identification_T0}, we obtain
\[
V^T(T-t,x,i) \ = \ V^{T'}(T'-t,x,i), \qquad \forall\,(t,x,i)\in[0,T]\times\R^d\times\I_m.
\]
This implies that the function $V$ given by \eqref{V=V^T} is well-defined. Moreover, $V$ is continuous and satisfies a growth condition as in \eqref{GrowthCond_T}. In addition, from the viscosity properties of $V^T$ it follows that $V$ is a viscosity solution to system \eqref{HJBintro} on $[0,T]\times\R^d\times\I_m$, for any $T>0$. From the arbitrariness of $T$, we have that $V$ is a viscosity solution to \eqref{HJBintro} on 
$[0,\infty)\times\R^d\times\I_m$.

\vspace{1mm}

\noindent\textbf{Step II.} \emph{Uniqueness.} Let $W\colon[0,\infty)\times\R^d\times\I_m\rightarrow\R$ be a continuous viscosity solution to \eqref{HJBintro} satisfying a growth condition as in \eqref{GrowthCond_T}. For any $T>0$ define the function $W^T\colon[0,T]\times\R^d\times\I_m\rightarrow\R$ as follows
\[
W^T(t,x,i) \ := \ W(T-t,x,i), \qquad \forall\,(t,x,i)\in[0,T]\times\R^d\times\I_m.
\]
Then $W^T$ is a continuous viscosity solution to \eqref{HJB_T} satisfying a polynomial growth condition as in \eqref{GrowthCond_V^T}. From comparison Theorem 4.1 in \cite{BCP14} it follows that $W^T\equiv V^T$, therefore $W\equiv V$.
\ep

\subsection{Long time asymptotics}

We present the following stochastic control representation bounds for every viscosity solution $(\lambda,\phi)$ to \eqref{HJB_ergodic}, with $\phi$ satisfying a polynomial growth condition.

\begin{Proposition} \label{P:repphi} 
Suppose that Assumption {\bf (H1)} holds and consider a viscosity solution $(\lambda,\phi)$ to \eqref{HJB_ergodic}. Then, for any $T>0$, $\phi$ is a viscosity supersolution $($resp. subsolution$)$ to the nonlinear backward parabolic system of variational inequalities in the unknown $\psi\colon[0,T]\times\R^d\times\I_m\rightarrow\R$$:$
\beq
\label{HJB_T_psi}
\begin{cases}
\min\Big\{- \dfrac{\partial\psi}{\partial t}(t,x,i) - \inf_{u\in U}\big[\Lc^{i,u}\psi(t,x,i) + f(x,i,u) - \lambda\big], \\
\psi(t,x,i) - \max_{j\neq i}\big[\psi(t,x,j) - c(x,i,j)\big]\Big\} \ = \ 0, &\!\!\!\!\!\!\!\!\!\!\!\!\!\!\!\!(t,x,i)\in[0,T)\times\R^d\times\I_m, \\
\psi(T,x,i) \ = \ \min_j \phi(x,j)\quad(\text{resp. }\max_j\phi(x,j)), &\!\!\!\!\!\!\!\!\!\!\!\!\!\!\!\!(x,i)\in\R^d\times\I_m.
\end{cases}
\enq
Suppose, in addition, that $\phi$ satisfies
\begin{equation}
\label{GrowthCond_phi}
|\phi(x,i)| \ \leq \ M_\phi(1 + |x|^{q_\phi}), \qquad \forall\,(x,i)\in\R^d\times\I_m,
\end{equation}
for some constants $M_\phi\geq0$ and $q_\phi\geq1$. Then, for any $T>0$, $\phi$ satisfies
\begin{equation}
\label{phi}
\sup_{\alpha\in\Ac_{t,t}}\inf_{\upsilon\in\Uc_{t,t}} \underline J_{(\lambda,\phi)}^T(t,x,i;\alpha,\upsilon) \ \leq \ \phi(x,i) \ \leq \ \sup_{\alpha\in\Ac_{t,t}}\inf_{\upsilon\in\Uc_{t,t}} \overline J_{(\lambda,\phi)}^T(t,x,i;\alpha,\upsilon),
\end{equation}
$\forall\,(t,x,i)\in[0,T]\times\R^d\times\I_m$, with
\beqs
\underline J_{(\lambda,\phi)}^T(t,x,i;\alpha,\upsilon) &:=& \E\bigg[\int_t^T \big(f(X_s^{t,x,i;\alpha,\upsilon},I_s^{t,x,i;\alpha,\upsilon},\upsilon_s) - \lambda\big)ds + \min_j\phi(X_T^{t,x,i;\alpha,\upsilon},j) \notag \\
& & \quad - \; \sum_{n\in\N} c(X_{\tau_n}^{t,x,i;\alpha,\upsilon},I_{\tau_n^-}^{t,x,i;\alpha,\upsilon},I_{\tau_n}^{t,x,i;\alpha,\upsilon}) 1_{\{\tau_n<T\}} \bigg], \\
\overline J_{(\lambda,\phi)}^T(t,x,i;\alpha,\upsilon) &:=& \E\bigg[\int_t^T \big(f(X_s^{t,x,i;\alpha,\upsilon},I_s^{t,x,i;\alpha,\upsilon},\upsilon_s) - \lambda\big)ds + \max_j\phi(X_T^{t,x,i;\alpha,\upsilon},j) \notag \\
& & \quad - \; \sum_{n\in\N} c(X_{\tau_n}^{t,x,i;\alpha,\upsilon},I_{\tau_n^-}^{t,x,i;\alpha,\upsilon},I_{\tau_n}^{t,x,i;\alpha,\upsilon}) 1_{\{\tau_n<T\}} \bigg].
\enqs
\end{Proposition}
\begin{Remark}
{\rm
Notice that from the viscosity properties of $\phi$ we also know that $\phi$ is a viscosity solution to the system with terminal condition $\phi$ itself:
\beqs
\begin{cases}
\min\Big\{- \dfrac{\partial\psi}{\partial t}(t,x,i) - \inf_{u\in U}\big[\Lc^{i,u}\psi(t,x,i) + f(x,i,u) - \lambda\big], \\
\psi(t,x,i) - \max_{j\neq i}\big[\psi(t,x,j) - c(x,i,j)\big]\Big\} \ = \ 0, &\!\!\!\!\!\!\!\!\!\!\!\!\!\!\!\!(t,x,i)\in[0,T)\times\R^d\times\I_m, \\
\psi(T,x,i) \ = \ \phi(x,i), &\!\!\!\!\!\!\!\!\!\!\!\!\!\!\!\!(x,i)\in\R^d\times\I_m.
\end{cases}
\enqs
Now, suppose that $\phi$ satisfies \eqref{GrowthCond_phi} and condition {\bf (H1)}(iii), i.e.,
\begin{equation}
\label{phi_terminal_condition}
\phi(x,i) \ \geq \ \max_{j\neq i}\big[\phi(x,j) - c(x,i,j)\big], \qquad \forall\,(x,i)\in\R^d\times\I_m.
\end{equation}
Then, from Corollary 4.1 in \cite{BCP14} it follows that $\phi$ admits the representation
\beq
\label{phi=J}
\phi(x,i) &=& \sup_{\alpha\in\Ac_{t,t}}\inf_{\upsilon\in\Uc_{t,t}} J_{(\lambda,\phi)}^T(t,x,i;\alpha,\upsilon), \qquad \forall\,(t,x,i)\in[0,T]\times\R^d\times\I_m,
\enq
with
\beqs
J_{(\lambda,\phi)}^T(t,x,i;\alpha,\upsilon) &:=& \E\bigg[\int_t^T \big(f(X_s^{t,x,i;\alpha,\upsilon},I_s^{t,x,i;\alpha,\upsilon},\upsilon_s) - \lambda\big)ds + \phi(X_T^{t,x,i;\alpha,\upsilon},I_T^{t,x,i;\alpha,\upsilon}) \notag \\
& & \quad - \; \sum_{n\in\N} c(X_{\tau_n}^{t,x,i;\alpha,\upsilon},I_{\tau_n^-}^{t,x,i;\alpha,\upsilon},I_{\tau_n}^{t,x,i;\alpha,\upsilon}) 1_{\{\tau_n<T\}} \bigg].
\enqs
However, since $(\lambda,\phi)$ is only a viscosity solution to \eqref{HJB_ergodic}, it is not obvious that \eqref{phi_terminal_condition} holds. For this reason, we introduce the two systems in \eqref{HJB_T_psi} with terminal conditions $\min_j\phi(x,j)$ and $\max_j\phi(x,j)$, which clearly satisfy condition \eqref{phi_terminal_condition}, since $c$ is nonnegative.
\ep
}
\end{Remark}
\textbf{Proof.} The fact that $\phi$ is a viscosity super/subsolution to \eqref{HJB_T_psi} follows obviously from the viscosity properties of $\phi$, since $\phi$ does not depend on time $t$.

From Corollary 4.1 in \cite{BCP14} we know that there exists a unique continuous and with polynomial growth viscosity solution $\underline\psi$ (resp. $\overline\psi$) to the system of variational inequalities \eqref{HJB_T_psi} with terminal condition $\min_j\phi(x,j)$ (resp. $\max_j\phi(x,j)$). From comparison Theorem 4.1 in \cite{BCP14} we have
\[
\underline\psi \ \leq \ \phi \ \leq \ \overline\psi.
\]
Using again Corollary 4.1 in \cite{BCP14}, we see that $\underline\psi$ and $\overline\psi$ admit the stochastic control representations
\[
\underline\psi(t,x,i) \ = \ \sup_{\alpha\in\Ac_{t,t}}\inf_{\upsilon\in\Uc_{t,t}} \underline J_{(\lambda,\phi)}^T(t,x,i;\alpha,\upsilon), \quad\; \overline\psi(t,x,i) \ = \ \sup_{\alpha\in\Ac_{t,t}}\inf_{\upsilon\in\Uc_{t,t}} \overline J_{(\lambda,\phi)}^T(t,x,i;\alpha,\upsilon),
\]
$\forall\,(t,x,i)\in[0,T]\times\R^d\times\I_m$, from which \eqref{phi} follows.
\ep

\vspace{3mm}

From the representation formula for $V$ and the bounds on a generic viscosity solution $(\lambda,\phi)$ to \eqref{HJB_ergodic}, we deduce the following result when $\phi$ satisfies a polynomial growth condition of second degree.

\begin{Proposition}
\label{P:V-->lambda}
Suppose that Assumptions {\bf (H1)}-{\bf (H2)} hold and consider a viscosity solution $(\lambda,\phi)$ to \eqref{HJB_ergodic}, with $\phi$ satisfying
\begin{equation}
\label{GrowthCond_phi2}
|\phi(x,i)| \ \leq \ M_\phi(1 + |x|^2), \qquad \forall\,(x,i)\in\R^d\times\I_m,
\end{equation}
for some constant $M_\phi\geq0$. Then
\[
\frac{V(T,x,i)}{T} \ \overset{T\rightarrow\infty}{\longrightarrow} \ \lambda, \qquad \forall\,(x,i)\in\R^d\times\I_m.
\]
In particular, $\lambda$ is uniquely determined for all viscosity solutions $(\lambda,\phi)$ to \eqref{HJB_ergodic}, with $\phi$ satisfying a polynomial growth condition as in \eqref{GrowthCond_phi2}.
\end{Proposition}
\begin{Remark}
{\rm
Notice that Proposition \ref{P:V-->lambda} gives the uniqueness of $\lambda$ and one of the convergence results claimed in Theorem \ref{T:Main}.
\ep
}
\end{Remark}
\textbf{Proof.}
From  \eqref{V=V^T} we have $V(T,x,i)=V^T(0,x,i)$. Therefore, using \eqref{V^T} we obtain
\[
V(T,x,i) \ = \ \sup_{\alpha\in\Ac_{0,0}}\inf_{\upsilon\in\Uc_{0,0}} J^T(0,x,i;\alpha,\upsilon).
\]
On the other hand, taking $t=0$ in \eqref{phi} we find
\[
\sup_{\alpha\in\Ac_{0,0}}\inf_{\upsilon\in\Uc_{0,0}} \underline J_{(\lambda,\phi)}^T(0,x,i;\alpha,\upsilon) \ \leq \ \phi(x,i) \ \leq \ \sup_{\alpha\in\Ac_{0,0}}\inf_{\upsilon\in\Uc_{0,0}} \overline J_{(\lambda,\phi)}^T(0,x,i;\alpha,\upsilon).
\]
Therefore
\[
|V(T,x,i) - \lambda T - \phi(x,i)| \leq \sup_{\alpha\in\Ac_{0,0}}\sup_{\upsilon\in\Uc_{0,0}} \E\big[\big|g(X_T^{0,x,i;\alpha,\upsilon},I_T^{0,x,i;\alpha,\upsilon})\big| + \max_j\big|\phi(X_T^{0,x,i;\alpha,\upsilon},j)\big|\big].
\]
From the growth condition on $g$ and $\phi$, and using estimate \eqref{EstimateX_diss}, it follows that there exists some positive constant $C$, independent of $T,x,i$, such that
\[
|V(T,x,i) - \lambda T - \phi(x,i)| \ \leq \ C(1+|x|^2),
\]
from which the claim follows.
\ep

\section{Asymptotic behavior of the elliptic system} \label{S:Elliptic} 

\setcounter{equation}{0} 
\setcounter{Theorem}{0} \setcounter{Proposition}{0}
\setcounter{Corollary}{0} \setcounter{Lemma}{0}
\setcounter{Definition}{0} \setcounter{Remark}{0}

In the present section, we study the elliptic system of variational inequalities \eqref{HJB_ellipticintro}, and in particular the asymptotic behavior of $V^\beta$ as $\beta$ tends to zero, which will entail the existence of a viscosity solution $(\lambda,\phi)$ to the ergodic system \eqref{HJB_ergodic}. Similarly to the parabolic case, we can derive a stochastic control representation for $V^\beta$ in terms of an infinite horizon robust feedback switching control problem. However, as we shall emphasize below, it is convenient to derive also another representation for $V^\beta$, known as dual representation formula, inspired by \cite{KP12} and based on randomization of the controls $\alpha$ and $\upsilon$.

\subsection{Motivation for the randomization approach}
\label{S:Motivation}

We begin presenting the representation of $V^\beta$ as value function of an infinite horizon robust feedback switching control, which reads as follows:
\beq
\label{V^beta}
V^\beta(x,i) &:=& \sup_{\alpha\in\Ac}\inf_{\upsilon\in\Uc} J^\beta(x,i;\alpha,\upsilon), \qquad \forall\,(x,i)\in\R^d\times\I_m,
\enq
with
\beqs
J^\beta(x,i;\alpha,\upsilon) &:=& \E\bigg[\int_0^\infty e^{-\beta t} f(X_t^{0,x,i;\alpha,\upsilon},I_t^{0,x,i;\alpha,\upsilon},\upsilon_t)dt \notag \\
& & \quad - \; \sum_{n\in\N} e^{-\beta\tau_n} c(X_{\tau_n}^{0,x,i;\alpha,\upsilon},I_{\tau_n^-}^{0,x,i;\alpha,\upsilon},I_{\tau_n}^{0,x,i;\alpha,\upsilon}) 1_{\{\tau_n<\infty\}} \bigg],
\enqs
where $\tau^n$ stands for $\tau^n(X_\cdot^{0,x,i;\alpha,\upsilon},I_{\cdot^-}^{0,x,i;\alpha,\upsilon})$, and the state processes $X^{0,x,i;\alpha,\upsilon},I^{0,x,i;\alpha,\upsilon}$ satisfy the first two equations in \eqref{SDE} with $t=0$ and for any $T\geq0$. Here, we do not pause on the technicalities, since the formulation of the stochastic control problem is analogous, with some obvious modifications, to that of the finite horizon case, already recalled in Section \ref{SubS:Robust}. We just notice that $\Ac$ denotes the family of all \emph{feedback switching controls} starting at time $0$ for the switcher. $\Ac$ is defined as $\Ac_{0,0}$ in Definition 2.1 of \cite{BCP14}, with $T$ replaced by $\infty$ and the property ``$\tau_n(y_n)=T$ for $n$ large enough along every adaptive sequence $(y_n)_n$'' is replaced by ``$\tau_n(y_n)\nearrow\infty$ as $n$ goes to infinity, for every adaptive sequence $(y_n)_n$''. On the other hand, $\Uc$ denotes the family of all \emph{open-loop controls} starting at time $0$ for nature. A generic element of $\Uc$ is an adapted process $\upsilon\colon[0,\infty)\times\Omega\rightarrow U$. It is then easy to see that, for every $T>0$, $(x,i)\in\R^d\times\I_m$, $\alpha\in\Ac$, $\upsilon\in\Uc$, Lemma \ref{L:X} still holds for $(X_s^{0,x,i;\alpha,\upsilon},I_s^{0,x,i;\alpha,\upsilon})_{0\leq s \leq T}$ and for the first two equations in \eqref{SDE} on $[0,T]$.

\begin{Proposition}
\label{P:VbetaStochasticPerron}
Let Assumptions {\bf (H1)} and {\bf (H2)} hold. Then, for any $\beta>0$, the function $V^\beta$ defined by \eqref{V^beta} is the unique continuous viscosity solution to system \eqref{HJB_ellipticintro} satisfying a linear growth condition
\begin{equation}\label{V^beta_linear_growth}
\sup_{(x,i)\in\R^d\times\I_m} \frac{|V^\beta(x,i)|}{1 + |x|} \ < \ \infty.
\end{equation}
\end{Proposition}
We do not report the proof of Proposition \ref{P:VbetaStochasticPerron} here (which provides statement (ii) of Theorem \ref{T:ExistHJB_T_beta}), since it can be done proceeding as in the finite horizon case, for which we refer to Corollary 4.1 in \cite{BCP14}. We just observe that the proof is based on the stochastic Perron method, which yields, as a by-product, the following dynamic programming principle:
\begin{align}
\label{V^betaDPP}
\hspace{-1mm}V^\beta(x,i) &= \sup_{\alpha\in\Ac}\inf_{\upsilon\in\Uc} \E\bigg[e^{-\beta T}V^\beta(X_T^{0,x,i;\alpha,\upsilon},I_T^{0,x,i;\alpha,\upsilon}) + \int_0^T e^{-\beta t} f(X_t^{0,x,i;\alpha,\upsilon},I_t^{0,x,i;\alpha,\upsilon},\upsilon_t)dt \notag \\
&\quad - \sum_{n\in\N} e^{-\beta\tau_n} c(X_{\tau_n}^{0,x,i;\alpha,\upsilon},I_{\tau_n^-}^{0,x,i;\alpha,\upsilon},I_{\tau_n}^{0,x,i;\alpha,\upsilon}) 1_{\{\tau_n\leq T\}} \bigg], \quad \forall\,(x,i)\in\R^d\times\I_m,
\end{align}
for all $T\geq0$. Identity \eqref{V^betaDPP} implies in particular the inequality
\begin{equation}
\label{VbetaInequality}
V^{\beta}(x,i) \ \geq \ \max_{j\neq i}\big[V^\beta(x,j) - c(x,i,j)\big], \qquad \forall\,(x,i)\in\R^d\times\I_m.
\end{equation}
As a matter of fact, fix $(x,i)\in\R^d\times\I_m$ and take $\bar\alpha=(\bar\tau_n,\bar\iota_n)_{n\in\N}$, with $\bar\tau_0\equiv0$, $\bar\tau_n\equiv+\infty$ for every $n\geq1$, and $\bar\iota_n\equiv j$ for every $n\in\N$. Then, from \eqref{V^betaDPP} we get
\begin{align}
V^\beta(x,i) \ &\geq \ \inf_{\upsilon\in\Uc} \E\bigg[e^{-\beta T}\big(V^\beta(X_T^{0,x,i;\bar\alpha,\upsilon},j) - V^\beta(x,j)\big) + \int_0^T e^{-\beta t} f(X_t^{0,x,i;\bar\alpha,\upsilon},j,\upsilon_t)dt\bigg] \notag \\
&\quad \ + e^{-\beta T}V^\beta(x,j) - c(x,i,j). \label{V-Vbis}
\end{align}
Since $V^\beta(\cdot,j)$ is a continuous function, for every $\eps>0$ there exists $\delta_\eps>0$ such that $|V^\beta(x,j)-V^\beta(x',j)|\leq\eps$, whenever $|x-x'|\leq\delta_\eps$. Therefore
\begin{align*}
&\sup_{\upsilon\in\Uc}\E\big[e^{-\beta T}\big|V^\beta(X_T^{0,x,i;\bar\alpha,\upsilon},j) - V^\beta(x,j)\big|\big] \\
&\leq \ e^{-\beta T}\eps + e^{-\beta T}\sup_{\upsilon\in\Uc}\E\big[\big|V^\beta(X_T^{0,x,i;\bar\alpha,\upsilon},j) - V^\beta(x,j)\big|1_{\{|X_T^{0,x,i;\bar\alpha,\upsilon} - x|>\delta_\eps\}}\big].
\end{align*}
Now, from \eqref{V^beta_linear_growth} we see that there exists some positive constant $c$ such that
\begin{align}
&\sup_{\upsilon\in\Uc}\E\big[e^{-\beta T}\big|V^\beta(X_T^{0,x,i;\bar\alpha,\upsilon},j) - V^\beta(x,j)\big|\big] \notag \\
&\leq \ e^{-\beta T}\eps + ce^{-\beta T}\sup_{\upsilon\in\Uc}\big\{\E\big[1 + |x|^2 + |X_T^{0,x,i;\bar\alpha,\upsilon}|^2\big]\P\big(|X_T^{0,x,i;\bar\alpha,\upsilon} - x|>\delta_\eps\big)\big\} \notag \\
&\leq \ e^{-\beta T}\eps + \frac{ce^{-\beta T}}{\delta_\eps^2}\sup_{\upsilon\in\Uc}\big\{\E\big[1 + |x|^2 + |X_T^{0,x,i;\bar\alpha,\upsilon}|^2\big]\E\big[|X_T^{0,x,i;\bar\alpha,\upsilon} - x|^2\big]\big\}. \label{V-V}
\end{align}
From estimate \eqref{EstimateX}, we see that $\sup_{\upsilon\in\Uc}\E[|X_T^{0,x,i;\bar\alpha,\upsilon}|^2]$ is finite. We also notice that the following standard estimate holds:
\[
\sup_{\upsilon\in\Uc}\E\big[|X_T^{0,x,i;\bar\alpha,\upsilon} - x|^2\big] \ \leq \ C'T(1 + |x|^2),
\]
for some positive constant $C'$. As a consequence, letting $T\downarrow0$ in \eqref{V-V}, we get
\[
\sup_{\upsilon\in\Uc}\E\big[e^{-\beta T}\big|V^\beta(X_T^{0,x,i;\bar\alpha,\upsilon},j) - V^\beta(x,j)\big|\big] \ \overset{T\downarrow0}{\longrightarrow} \ 0.
\]
Similarly
\[
\sup_{\upsilon\in\Uc}\E\bigg[\int_0^T e^{-\beta t} \big|f(X_t^{0,x,i;\bar\alpha,\upsilon},j,\upsilon_t)\big|dt\bigg] \ \overset{T\downarrow0}{\longrightarrow} \ 0.
\]
In conclusion, letting $T\downarrow0$ in \eqref{V-Vbis} we obtain
\[
V^\beta(x,i) \ \geq \ V^\beta(x,j) - c(x,i,j).
\]
From the arbitrariness of $j$, we deduce that \eqref{VbetaInequality} holds.

Inequality \eqref{VbetaInequality} is the only result, derived from the stochastic control representation \eqref{V^beta}, that we shall exploit. Instead, for all the other results, we shall use the dual representation formula for $V^\beta$. To understand why, we begin noting that, because of the feedback form of $\alpha$, from \eqref{V^beta} it is not clear how to prove some properties of $V^\beta$, which are crucial to perform the asymptotic analysis. As an example, for the proof of the uniform Lipschitz property of $V^\beta$ with respect to $x$ we need an estimate of the following type (see estimate \eqref{EstimateX-X'}): for any $t\geq0$, there exists a constant $L_t\geq0$, depending only on $t$, such that
\begin{equation}
\label{Estimate_feedback}
\sup_{i\in\I_m,\alpha\in\Ac,\upsilon\in\Uc}\E\big[|X_t^{0,x,i;\alpha,\upsilon} - X_t^{0,x',i;\alpha,\upsilon}|\big] \ \leq \ L_t|x-x'|, \qquad\forall\,x,x'\in\R^d.
\end{equation}
It is however not clear how to prove estimate \eqref{Estimate_feedback} from representation \eqref{V^beta}, since the evolution of the process $I$ (which influence the dynamics of $X$ through the coefficients $b$ and $\sigma$) depends on the starting point, $x$ or $x'$, due to the feedback form of $\alpha$.

As a possible solution to the issue raised above, one could consider the Elliott-Kalton version of the infinite horizon robust switching control problem (which can be formulated in a similar way to the finite horizon case, for which we refer to Section 4.2 in \cite{BCP14}). More precisely, suppose we are able to prove that $V^\beta$ is the value function of the Elliott-Kalton version of our control problem. Since non-anticipative strategies for the switcher take values in the class of all switching controls (not necessarily of feedback type), estimate \eqref{Estimate_feedback} follows easily from that representation. The drawback of this approach is the proof that $V^\beta$ is indeed the value function of the Elliott-Kalton formulation, which in general relies on delicate measurability issues arising in the proof of the dynamic programming principle, firstly faced in the seminal paper \cite{fleming_souganidis}.

An alternative to the Elliott-Kalton representation of $V^\beta$ is a dual representation formula (in the sense of \cite{KP12}, Section 2.3) based on randomization of the controls $\alpha$ and $\upsilon$ (see Remarks \ref{R:Randomization} and \ref{R:DualProblem} for some insights on the dual representation formula). This latter turns out to be easier to derive, since it avoids the proof of a dynamic programming principle, but it still allows to prove estimates for the process $X$ as \eqref{Estimate_feedback} above. For this reason, our aim is now to prove that $V^\beta$ admits a dual representation formula, which can be deduced starting from an opportune penalized elliptic system of variational inequalities that we now introduce.

\subsection{Penalized elliptic system of variational inequalities}

For any $n\in\N$ and $\beta>0$, consider the following penalized system of variational inequalities in the unknown $V^{\beta,n}\colon\R^d\times\I_m\rightarrow\R$:
\begin{align}
\label{HJB_penalized}
\beta V^{\beta,n}(x,i) - \inf_{u\in U}\big[\Lc^{i,u}V^{\beta,n}(x,i) + f(x,i,u)\big]\;\;\,& \\
- n\sum_{j=1}^m \big[V^{\beta,n}(x,j) - V^{\beta,n}(x,i)
- c(x,i,j)\big]^+& \ = \ 0, \notag
\end{align}
for any $(x,i)\in\R^d\times\I_m$, where $h^+=\max(h,0)$ denotes the positive part of the function $h$. Our aim is to prove that there exists a continuous viscosity solution $V^{\beta,n}$ to system \eqref{HJB_penalized}, converging to $V^\beta$ as $n$ tends to infinity, such that $V^{\beta,n}$ satisfies a linear growth condition and other opportune estimates uniformly in $n$, so that they still hold for $V^\beta$ letting $n\rightarrow\infty$. We postpone the proof of the convergence of $V^{\beta,n}$ towards $V^\beta$, and we begin focusing on the proof of the estimates for $V^{\beta,n}$, for which we adopt a probabilistic approach. In particular, inspired by the results in \cite{KP12}, we derive Feynman-Kac and dual representation formulae of $V^{\beta,n}$ by means of a backward stochastic differential equation with partially nonnegative jumps on infinite horizon.

\subsubsection{Feynman-Kac and dual representation formulae of $V^{\beta,n}$}

We begin introducing some notations (to simplify the presentation, we redefine certain symbols as $\Omega$ or $\F$, already used in Section \ref{S:Parabolic}, since no confusion should arise).

\vspace{1mm}

Let $(\Omega,\Fc,\P)$ be a complete probability space, $W$ $=$ $(W_t)_{t\geq 0}$ a $d$-dimensional Brownian motion, $\pi$ a Poisson random measure on $\R_+\times \I_m$, $\mu$ a Poisson random measure on $\R_+\times U$, such that $W$, $\pi$, and $\mu$ are independent. We assume that $\pi$ has intensity measure $\vartheta_\pi(di)dt$, where $\vartheta_\pi(di)=\sum_{j=1}^m \delta_j(di)$ and $\delta_j$ denotes the Dirac delta at $j\in\I_m$. We also suppose that $\mu$ has intensity measure $\vartheta_\mu(du)dt$, where $\vartheta_\mu$ is a finite measure on $(U,\Bc(U))$ ($\Bc(U)$ denotes the Borel $\sigma$-field on $U$) such that:
\begin{itemize}
\item[(i)] The support of $\vartheta_\mu$ is the whole set $\mathring U$, namely $\vartheta_\mu(\Oc\cap\mathring U)>0$ for any open subset $\Oc$ of $\R^q$ satisfying $\Oc\cap\mathring U\neq\emptyset$.
\item[(ii)] The boundary $\partial U=U\setminus\mathring U$ of $U$ is $\vartheta_\mu$-negligible: $\vartheta_\mu(\partial U)=0$. 
\end{itemize}
We denote by $\tilde\pi(dt,di)$ $=$ $\pi(dt,di)-\vartheta_\pi(di)dt$ and $\tilde\mu(dt,du)$ $=$ $\mu(dt,du)-\vartheta_\mu(du)dt$ the compensated martingale measures
associated with $\pi$ and $\mu$, respectively. We also set $\F$ $=$ $(\Fc_t)_{t\geq 0}$ the completion of the natural filtration generated by $W$, $\pi$, and $\mu$, and we denote by $\Pc$ the $\sigma$-field of $\F$-predictable subsets of $\Omega\times [0,\infty)$.

\paragraph{Forward SDE and randomization.} For any $(x,i,u)\in\R^d\times\I_m\times U$, we introduce the forward jump-diffusion Markov process $(X,I,\Gamma)$, which evolves on $[0,\infty)$ according to the following system of stochastic differential equations:
\begin{equation}
\label{FSDE}
\begin{cases}
\displaystyle X_t = x + \int_0^t b(X_s,I_s,\Gamma_s) ds + \int_0^t \sigma(X_s,I_s,\Gamma_s) dW_s, \\
\displaystyle I_t = i + \int_0^t \int_{\I_m} (j - I_{s^-}) \pi(ds,d\,\!j), \\
\displaystyle\Gamma_t = u + \int_0^t \int_U (u' - \Gamma_{s^-}) \mu(ds,du'),
\end{cases}
\end{equation}
for all $t\geq0$.
\begin{Remark}
\label{R:Randomization}
{\rm
The process $(X,I,\Gamma)$ satisfying system \eqref{FSDE} will be the forward and driving process in the Feynman-Kac and dual representation formulae of $V^{\beta,n}$. 
The expression of \eqref{FSDE} is derived following \cite{KP12} and it is inspired by the stochastic control representation of $V^\beta$.  As stated in Proposition \ref{P:VbetaStochasticPerron}, $V^\beta$ is related to an infinite horizon robust switching control problem, whose state process evolves according to the first two equations in \eqref{SDE} (with $t=0$ and for any $T\in[0,\infty)$). Then, \eqref{FSDE} is obtained from \eqref{SDE} by means of a \emph{randomization} of the switching and open-loop controls, which is fulfilled introducing the pure jump Markov processes $I$ and $\Gamma$ driven by independent Poisson random measures. Hence,  compared to \reff{SDE}, the pair $(X,I)$ in \eqref{FSDE}   is an uncontrolled process. 
\ep
}
\end{Remark}

It is well-known that under Assumption {\bf (H1)}, for any $(x,i,u)\in\R^d\times\I_m\times U$, there exists a unique solution $(X^{x,i,u},I^i,\Gamma^u)=(X_t^{x,i,u},I_t^i,\Gamma_t^u)_{t\geq0}$ to system \eqref{FSDE}, and the following standard estimate holds: for any $T\geq0$ and $p\geq2$, there exists some positive constant $C_{p,T}$ such that
\[
\E\Big[\sup_{0\leq t\leq T}|X_t^{x,i,u}|^p\Big] \ \leq \ C_{p,T}\big(1 + |x|^p\big).
\]
We conclude this paragraph with two important estimates resulting from the dissipativity condition {\bf (H2)}, which call on a family of probability measures we are going to define. For any $n\in\N\backslash\{0\}$ and $k\in\N$, let $\Xi_n$ be the set of $\Pc\otimes\mathscr P(\I_m)$-measurable\footnote{$\mathscr P(\I_m)$ denotes the power set of $\I_m$.} maps valued in $(0,n]$, $\Vc_k$ be the set of $\Pc\otimes\Bc(U)$-measurable maps valued in $[1,k+1]$, and denote by $\Xi=\cup_{n\in\N\backslash\{0\}}\Xi_n$, $\Vc=\cup_{k\in\N}\Vc_k$. We consider for $\xi\in\Xi$, $\nu\in\Vc$, the probability measure $\P^{\xi,\nu}$ equivalent to $\P$ on $(\Omega,\Fc_T)$, for any $T>0$, with Radon-Nikodym density:
\begin{equation}
\label{P^xi,nu}
\frac{d\P^{\xi,\nu}}{d\P} \Big|_{\Fc_T} \!\! = \ \zeta_T^{\xi,\nu} \ := \ \Ec_T\Big( \int_0^. \int_{\I_m} \! (\xi_t(j) -1) \tilde\pi(dt,d\,\!j) \Big) \cdot\, \Ec_T\Big( \int_0^. \int_U \! (\nu_t(u') - 1) \tilde\mu(dt,du') \Big),
\end{equation}
where $\Ec(\cdot)$ denotes the Dol\'eans-Dade exponential local martingale. Proceeding as in Lemma 2.4 of \cite{KP12}, we see that $\zeta^{\xi,\nu}$ is a ``true" $\P$-martingale (hence defining a probability measure $\P^{\xi,\nu}$ through \eqref{P^xi,nu}) since $\xi,\nu$ are essentially bounded and $\vartheta_\pi,\vartheta_\mu$ are finite measures on $\I_m$ and $U$, with $\zeta_T^{\xi,\nu}\in{\bf L^2}(\Omega,\Fc_T,\P)$ for any $T>0$. By Girsanov's theorem, we recall that $W$ remains a Brownian motion under $\P^{\xi,\nu}$, and the effect of the probability measure $\P^{\xi,\nu}$ is to change the compensator $\vartheta_\pi(d\,\!i)dt$ of $\pi$ under $\P$ to $\xi_t(i)\vartheta_\pi(d\,\!i)dt$ under
$\P^{\xi,\nu}$, and the compensator $\vartheta_\mu(du)dt$ of $\mu$ under $\P$ to $\nu_t(u)\vartheta_\mu(du)dt$ under  $\P^{\xi,\nu}$. We denote by $\tilde\pi^\xi(dt,d\,\!i)$ $=$ $\pi(dt,d\,\!i)-\xi_t(i)\vartheta_\pi(d\,\!i)dt$ and $\tilde\mu^\nu(dt,du)$ $=$ $\mu(dt,du)-\nu_t(u)\vartheta_\mu(du)dt$ the compensated martingale measures associated with $\pi$ and $\mu$ under $\P^{\xi,\nu}$. Finally, $\E^{\xi,\nu}$ denotes the expectation with respect to $\P^{\xi,\nu}$.

\begin{Lemma}
Let Assumptions {\bf (H1)} and {\bf (H2)} hold.
\begin{itemize}
\item[\textup{(i)}] For all $t\geq0$, $x,x'\in\R^d$, $i\in\I_m$, $u\in U$,
\beq
\label{EstimateX-X'}
\sup_{\xi\in\Xi,\,\nu\in\Vc}\E^{\xi,\nu}\big[|X_t^{x,i,u} - X_t^{x',i,u}|^2\big] &\leq & e^{-2\gamma t} |x - x'|^2.
\enq
\item[\textup{(ii)}] There exists a positive constant $\bar C$, depending only on $M_{b,\sigma}=\sup_{(i,u)\in\I_m\times U}(|b(0,i,u)|$ $+$ $\|\sigma(0,i,u)\|)$, $L_1$, and $\gamma$, such that
\beq
\label{EstimateX_diss_xi,nu}
\sup_{t\geq0,\,\xi\in\Xi,\,\nu\in\Vc}\E^{\xi,\nu}\big[|X_t^{x,i,u}|^2\big] &\leq & \bar C(1+|x|^2),
\enq
for any $(x,i,u)\in\R^d\times\I_m\times U$.
\end{itemize}
\end{Lemma}
\textbf{Proof.} (i) \emph{Estimate \eqref{EstimateX-X'}.} The proof of estimate \eqref{EstimateX-X'} can be done proceeding as in Lemma 2.1(ii) in \cite{CFP14}. Let us just give an idea of the proof. Firstly, we apply It\^o formula to $|X_s^{x,i,u} - X_s^{x',i,u}|^2$ between $s=0$ and $s=t\geq0$, then we take the expectation $\E^{\xi,\nu}$ with respect to $\P^{\xi,\nu}$, and finally we use the dissipativity condition \eqref{dissipative}. In conclusion, we end up with
\[
\E^{\xi,\nu}\big[|X_t^{x,i,u} - X_t^{x',i,u}|^2\big] \ \leq \ |x -x'|^2 - 2\gamma \int_0^t \E^{\xi,\nu}\big[|X_s^{x,i,u} - X_s^{x',i,u}|^2\big] ds.
\]
Then, the claim follows from Gronwall's inequality.

\vspace{1mm}

\noindent (ii) \emph{Estimate \eqref{EstimateX_diss_xi,nu}.} The result can be proved proceeding as in Lemma \ref{L:X_diss}.
\ep

\paragraph{Backward SDE with partially nonnegative jumps on infinite horizon.} For any $T\in[0,\infty)$, we denote by $\Pc_T$ the $\sigma$-field of $\F$-predictable subsets of $\Omega\times [0,T]$ and we introduce the following spaces of random maps:

\begin{itemize}
\item ${\bf S_T^2}$ the set of real-valued c\`adl\`ag $\F$-adapted processes $Y$ $=$
$(Y_t)_{0\leq t\leq T}$ satisfying
\[
\|Y\|_{_{{\bf S_T^2}}}^2 \ := \ \E\Big[ \sup_{0\leq t\leq T} |Y_t|^2 \Big] \ < \ \infty.
\]
We also denote ${\bf S_{\textup{loc}}^2}:=\cap_{T>0}{\bf S_T^2}$.
\item ${\bf L_T^p(W)}$, $p\geq1$, the set of $\R^d$-valued $\Pc_T$-measurable processes
$Z=(Z_t)_{0\leq t\leq T}$ satisfying
\[
\|Z\|_{_{\bf L_T^p(W)}}^p \ := \ \E\bigg[\bigg(\int_0^T |Z_t|^2 dt\bigg)^{\frac{p}{2}}\bigg] \ < \ \infty.
\]
We also denote ${\bf L_{\textup{loc}}^p(W)}:=\cap_{T>0}{\bf L_T^p(W)}$.
\item ${\bf L_T^p(\tilde\pi)}$, $p\geq1$, the set of
$\Pc_T\otimes\mathscr{P}(\I_m)$-measurable maps $L\colon\Omega\times[0,T]\times\I_m\rightarrow\R$ satisfying
\[
\|L\|_{_{{\bf L_T^p(\tilde\pi)}}}^p \ := \ \E\bigg[\int_0^T \bigg(\sum_{i=1}^m|L_t(i)|^2\bigg)^{\frac{p}{2}}dt\bigg] \ < \ \infty.
\]
We also denote ${\bf L_{\textup{loc}}^p(\tilde\pi)}:=\cap_{T>0}{\bf L_T^p(\tilde\pi)}$.
\item ${\bf L_T^p(\tilde\mu)}$ the set of
$\Pc_T\otimes\Bc(U)$-measurable maps $R\colon\Omega\times[0,T]\times U\rightarrow\R$ satisfying
\[
\|R\|_{_{{\bf L_T^p(\tilde\mu)}}}^p \ := \ \E\bigg[\int_0^T\bigg(\int_U |R_t(u)|^2 \vartheta_\mu(du)\bigg)^{\frac{p}{2}}dt\bigg] \ < \ \infty.
\]
We also denote ${\bf L_{\textup{loc}}^p(\tilde\mu)}:=\cap_{T>0}{\bf L_T^p(\tilde\mu)}$.
\item ${\bf K_T^2}$ the set of nondecreasing predictable processes $K$ $=$ $(K_t)_{0\leq t\leq T}$ $\in$  ${\bf S_T^2}$ with $K_0$ $=$ $0$, so that
\[
\|K\|_{_{{\bf S_T^2}}}^2 \ = \ \E|K_T|^2.
\]
We also denote ${\bf K_{\textup{loc}}^2}:=\cap_{T>0}{\bf K_T^2}$.
\end{itemize}
Let us now consider, for any $\beta>0$, $n\in\N$, and $(x,i,u)\in\R^d\times\I_m\times U$, the following backward stochastic differential equation with partially nonnegative jumps on infinite horizon:
\beq
\label{BSDE}
Y_t^{\beta,n} &=& Y_T^{\beta,n} - \beta\int_t^T Y_s^{\beta,n}ds + \int_t^T f(X_s^{x,i,u},I_s^i,\Gamma_s^u) ds - \sum_{j=1}^m \int_t^T L_s^{\beta,n}(j) ds \notag \\
&& + \; n\sum_{j=1}^m\int_t^T \big[L_s^{\beta,n}(j) - c(X_s^{x,i,u},I_{s^-}^i,j)\big]^+ ds - \big(K_T^{\beta,n} - K_t^{\beta,n}\big)  \\
&& - \; \int_t^T Z_s^{\beta,n} dW_s - \int_t^T\int_{\I_m} L_s^{\beta,n}(j) \tilde\pi(ds,d\,\!j) - \int_t^T\int_U R_s^{\beta,n}(u') \tilde\mu(ds,du'), \notag
\enq
for any $0\leq t\leq T$, $T\in[0,\infty)$, and
\beq
\label{BSDE_Constraint}
R_t^{\beta,n}(u') &\geq &0, \qquad d\P\otimes dt\otimes\vartheta_\mu(du')\text{-a.e.}
\enq
Our aim is to prove that the penalized elliptic system \eqref{HJB_penalized} is related to the \emph{maximal solution} to the above BSDE with partially nonnegative jumps \eqref{BSDE}-\eqref{BSDE_Constraint}, that we now define.

\begin{Definition}
\label{D:Maximal}
For any $\beta>0$, $n\in\N$, and $(x,i,u)\in\R^d\times\I_m\times U$, we say that a quintuple $(Y^{\beta,n,x,i,u},Z^{\beta,n,x,i,u},L^{\beta,n,x,i,u},R^{\beta,n,x,i,u},K^{\beta,n,x,i,u})\in{\bf S_{\textup{loc}}^2}\times{\bf L_{\textup{loc}}^2(W)}\times{\bf L_{\textup{loc}}^2(\tilde\pi)}\times{\bf L_{\textup{loc}}^2(\tilde\mu)}\times{\bf K_{\textup{loc}}^2}$ is a maximal solution to the BSDE with partially nonnegative jumps on infinite horizon \eqref{BSDE}-\eqref{BSDE_Constraint} if:
\begin{itemize}
\item[\textup{(i)}] $|Y_t^{\beta,n,x,i,u}|\leq C(1+|X_t^{x,i,u}|)$, for all $t\geq0$ and for some positive constant $C$.
\item[\textup{(ii)}] For any other solution $(\underline Y^{\beta,n},\underline Z^{\beta,n},\underline L^{\beta,n},\underline R^{\beta,n},\underline K^{\beta,n})\in{\bf S_{\textup{loc}}^2}\times{\bf L_{\textup{loc}}^2(W)}\times{\bf L_{\textup{loc}}^2(\tilde\pi)}\times{\bf L_{\textup{loc}}^2(\tilde\mu)}\times{\bf K_{\textup{loc}}^2}$ to \eqref{BSDE}-\eqref{BSDE_Constraint} satisfying $|\underline Y_t^{\beta,n}|\leq\underline C(1+|X_t^{x,i,u}|)$, $\forall\,t\geq0$ and for some positive constant $\underline C$ (possibly depending on $\beta,n,x,i,u$), we have
\[
Y_t^{\beta,n,x,i,u} \ \geq \ \underline Y_t^{\beta,n}, \qquad \P\text{-a.s., for all }t\geq0.
\] 
\end{itemize}
\end{Definition}

\begin{Proposition}
\label{T:Feynman-Kac}
Let Assumptions {\bf (H1)}, {\bf (H2)}, and {\bf (H$U$)} hold. Then, for any $\beta>0$ and $n\in\N$, we have$:$
\begin{itemize}
\item[\textup{(i)}] For any $(x,i,u)\in\R^d\times\I_m\times U$, there exists a unique maximal solution to \eqref{BSDE}-\eqref{BSDE_Constraint} denoted by $(Y^{\beta,n,x,i,u},Z^{\beta,n,x,i,u},L^{\beta,n,x,i,u},R^{\beta,n,x,i,u},K^{\beta,n,x,i,u})$.
\item[\textup{(ii)}] Given $(x,i)\in\R^d\times\I_m$, for any $u\in\mathring U$ the random variable $Y_0^{\beta,n,x,i,u}$ is equal $\P$-a.s. to a constant independent of $u\in\mathring U$. Moreover, the function $V^{\beta,n}\colon\R^d\times\I_m\rightarrow\R$ given by the \textbf{Feynman-Kac formula}:
\[
V^{\beta,n}(x,i) \ = \ Y_0^{\beta,n,x,i,u}, \qquad \forall\,(x,i)\in\R^d\times\I_m,
\]
for any $u\in\mathring U$, is a continuous viscosity solution to system \eqref{HJB_penalized} and satisfies the linear growth condition
\[
\sup_{(x,i)\in\R^d\times\I_m}\frac{|V^{\beta,n}(x,i)|}{1 + |x|} \ < \ \infty.
\]
\end{itemize}
\end{Proposition}
\textbf{Proof}
See Appendix, in particular Proposition \ref{P:PenalizedBSDE} for (i) and Proposition \ref{T:Feynman-Kac_App} for (ii).
\ep

\begin{Corollary}
\label{C:DualFormula}
Let Assumptions {\bf (H1)}, {\bf (H2)}, and {\bf (H$U$)} hold. Then, for any $\beta>0$ and $n\in\N\backslash\{0\}$, the following \textbf{dual representation formula} for $V^{\beta,n}$ holds:
\begin{align}
\label{DualFormulaV^beta,n}
V^{\beta,n}(x,i) \ = \ \sup_{\xi\in\Xi_n}\inf_{\nu\in\Vc}\E^{\xi,\nu}\bigg[&\int_0^\infty e^{-\beta t}f(X_t^{x,i,u},I_t^i,\Gamma_t^u) dt \\
& - \int_0^\infty\int_{\I_m} e^{-\beta t} c(X_t^{x,i,u},I_{t^-}^i,j) \pi(dt,d\,\!j) \bigg], \;\; \forall\,(x,i)\in\R^d\times\I_m, \notag
\end{align}
for any $u\in\mathring U$.
\end{Corollary}
\textbf{Proof.}
As a preliminary step, the present proof involves the doubly indexed penalized BSDE \eqref{BSDEnk} introduced in the Appendix. Indeed, we firstly derive a dual representation formula for $V^{\beta,n,k}(x,i,u)=Y_0^{\beta,n,k,x,i,u}$, then we pass to the limit as $k\rightarrow\infty$, using that $Y_0^{\beta,n,k,x,i,u}\searrow Y_0^{\beta,n,x,i,u}=V^{\beta,n}(x,i)$ $\P$-a.s., to derive \eqref{DualFormulaV^beta,n}.

\vspace{1mm}

\noindent\textbf{Step I.} \emph{Dual representation formula for $V^{\beta,n,k}$.} Fix $n\in\N\backslash\{0\}$, $k\in\N$, and $(x,i,u)\in\R^d\times\I_m\times U$. Our aim is to prove the following dual representation formula
\begin{align}
\label{DualFormulaV^beta,n,k}
V^{\beta,n,k}(x,i,u) \ = \ \sup_{\xi\in\Xi_n}\inf_{\nu\in\Vc_k}\E^{\xi,\nu}\bigg[&\int_0^\infty e^{-\beta t}f(X_t^{x,i,u},I_t^i,\Gamma_t^u) dt \\
&- \int_0^\infty\int_{\I_m} e^{-\beta t} c(X_t^{x,i,u},I_{t^-}^i,j) \pi(dt,d\,\!j) \bigg]. \notag
\end{align}
To this end, take $\xi\in\Xi_n$ and $\nu\in\Vc_k$. Then, given $T\in[0,\infty)$, we add and subtract to \eqref{BSDEnk}, with $t=0$, the two terms (we adopt the simplified notation $Y^{\beta,n,k}$ for $Y^{\beta,n,k,x,i,u}$, and similarly for the other components)
\[
\int_0^T\int_{\I_m} L_t^{\beta,n,k}(j) \xi_t(j)\vartheta_\pi(d\,\!j)dt, \qquad \int_0^T\int_U R_t^{\beta,n,k}(u') \nu_t(u')\vartheta_\mu(du')dt.
\]
From the integrability conditions on $L^{\beta,n,k}$ (resp. $R^{\beta,n,k}$), it follows from Proposition II.1.28 in \cite{jacod_shiryaev} that the stochastic integral in \eqref{BSDEnk} with respect to $\tilde\pi$ (resp. $\tilde\mu$) can be written as the difference between the integral with respect to $\pi$ (resp. $\mu$) and that with respect to its $\P$-compensator $\vartheta_\pi(d\,\!j)dt$ (resp. $\vartheta_\mu(du')dt$). The same remark applies to the stochastic integral of $L^{\beta,n,k}$ (resp. $R^{\beta,n,k}$) with respect to $\tilde\pi^\xi$ (resp. $\tilde\mu^\nu$). Then, rearranging the terms in \eqref{BSDEnk}, we end up with (recall that $\vartheta_\pi(di)=\sum_{j=1}^m\delta_j(di)$)
\beq
\label{BSDEnk_xi_nu}
Y_0^{\beta,n,k} &=& Y_T^{\beta,n,k} - \beta\int_0^T Y_t^{\beta,n,k}dt + \int_0^T f(X_t^{x,i,u},I_t^i,\Gamma_t^u) dt \\
&& + \; n\sum_{j=1}^m\int_0^T \! \big[L_t^{\beta,n,k}(j) - c(X_t^{x,i,u},I_{t^-}^i,j)\big]^+ dt - k\int_0^T \! \int_U \big[R_t^{\beta,n,k}(u')\big]^- \vartheta_\mu(du')dt \notag \\
&& - \; \int_0^T Z_t^{\beta,n,k} dW_t - \int_0^T\int_{\I_m} L_t^{\beta,n,k}(j) \tilde\pi^\xi(dt,d\,\!j) - \int_0^T\int_U R_t^{\beta,n,k}(u') \tilde\mu^\nu(dt,du') \notag \\
&& - \; \sum_{j=1}^m \int_0^T L_t^{\beta,n,k}(j)\xi_t(j) dt - \int_0^T\int_U R_t^{\beta,n,k}(u')(\nu_t(u') - 1)\vartheta_\mu(du')dt. \notag
\enq
Now, we apply It\^o formula to $e^{-\beta t}Y_t^{\beta,n,k}$ between $0$ and $T$, afterwards we add and subtract the term
\[
\sum_{j=1}^m\int_0^T e^{-\beta t} c(X_t^{x,i,u},I_{t^-}^i,j) \xi_t(j) dt.
\]
Therefore, from \eqref{BSDEnk_xi_nu} we obtain
\begin{align}
\label{BSDEnk_xi_nu2}
&Y_0^{\beta,n,k} \ = \ e^{-\beta T} Y_T^{\beta,n,k} + \int_0^T e^{-\beta t} f(X_t^{x,i,u},I_t^i,\Gamma_t^u) dt - \sum_{j=1}^m\int_0^T e^{-\beta t} c(X_t^{x,i,u},I_{t^-}^i,j) \xi_t(j) dt \notag \\
& + \sum_{j=1}^m\int_0^T e^{-\beta t}\big\{n\big[L_t^{\beta,n,k}(j) - c(X_t^{x,i,u},I_{t^-}^i,j)\big]^+ - \xi_t(j)\big[L_t^{\beta,n,k}(j) - c(X_t^{x,i,u},I_{t^-}^i,j)\big]\big\} dt \notag \\
& - \int_0^T\int_U e^{-\beta t}\big\{k\big[R_t^{\beta,n,k}(u')\big]^- + (\nu_t(u') - 1)R_t^{\beta,n,k}(u')\big\} \vartheta_\mu(du')dt \\
& - \int_0^T e^{-\beta t} Z_t^{\beta,n,k} dW_t - \int_0^T\int_{\I_m} e^{-\beta t} L_t^{\beta,n,k}(j) \tilde\pi^\xi(dt,d\,\!j) - \int_0^T\int_U e^{-\beta t} R_t^{\beta,n,k}(u') \tilde\mu^\nu(dt,du'). \notag
\end{align}
Reasoning as in Lemma 2.5 in \cite{KP12}, we can prove that the three stochastic integrals appearing in \eqref{BSDEnk_xi_nu2}, which are $\P^{\xi,\nu}$-local martingales, are indeed true $\P^{\xi,\nu}$-martingales. Therefore, taking the expectation $\E^{\xi,\nu}$ with respect to $\P^{\xi,\nu}$ in \eqref{BSDEnk_xi_nu2}, we find
\begin{align}
\label{BSDEnk_xi_nu3}
&Y_0^{\beta,n,k} \ = \ G_T(\xi,\nu) \\
& + \sum_{j=1}^m\int_0^T \! e^{-\beta t}\E^{\xi,\nu}\big\{n\big[L_t^{\beta,n,k}(j) - c(X_t^{x,i,u},I_{t^-}^i,j)\big]^+ - \xi_t(j)\big[L_t^{\beta,n,k}(j) - c(X_t^{x,i,u},I_{t^-}^i,j)\big]\big\} dt \notag \\
& - \int_0^T\int_U e^{-\beta t}\E^{\xi,\nu}\big\{k\big[R_t^{\beta,n,k}(u')\big]^- + (\nu_t(u') - 1)R_t^{\beta,n,k}(u')\big\} \vartheta_\mu(du')dt, \notag
\end{align}
where
\begin{align*}
G_T(\xi,\nu) \ &= \ \E^{\xi,\nu}\bigg[e^{-\beta T} Y_T^{\beta,n,k} + \int_0^T e^{-\beta t} f(X_t^{x,i,u},I_t^i,\Gamma_t^u) dt \\
&\quad \ - \sum_{j=1}^m\int_0^T e^{-\beta t} c(X_t^{x,i,u},I_{t^-}^i,j) \xi_t(j) dt\bigg].
\end{align*}
Let us prove the following identities
\begin{equation}
\label{DualFormula_T}
Y_0^{\beta,n,k} \ = \ \sup_{\xi\in\Xi_n}\inf_{\nu\in\Vc_k} G_T(\xi,\nu) \ = \ \inf_{\nu\in\Vc_k}\sup_{\xi\in\Xi_n} G_T(\xi,\nu).
\end{equation}
We begin noting that, for any $\xi\in\Xi_n$ and $\nu\in\Vc_k$, we have
\begin{align}
n\big[L_t^{\beta,n,k}(j) - c(X_t^{x,i,u},I_{t^-}^i,j)\big]^+ - \xi_t(j)\big[L_t^{\beta,n,k}(j) - c(X_t^{x,i,u},I_{t^-}^i,j)\big] \ &\geq \ 0, \label{Ineq_xi} \\
k\big[R_t^{\beta,n,k}(u')\big]^- + (\nu_t(u') - 1)R_t^{\beta,n,k}(u') \ &\geq \ 0. \label{Ineq_nu}
\end{align}
Now, for every $\eps>0$, define $\xi^\eps\in\Xi_n$ as follows
\[
\xi_t^\eps(j) \ = \
\begin{cases}
n, &\quad L_t^{\beta,n,k}(j) - c(X_t^{x,i,u},I_{t^-}^i,j)\geq0, \\
\eps,  &\quad -1\leq L_t^{\beta,n,k}(j) - c(X_t^{x,i,u},I_{t^-}^i,j)<0, \\
-\dfrac{\eps}{L_t^{\beta,n,k}(j) - c(X_t^{x,i,u},I_{t^-}^i,j)}, &\quad L_t^{\beta,n,k}(j) - c(X_t^{x,i,u},I_{t^-}^i,j)<-1.
\end{cases}
\]
On the other hand, let $\nu^*\in\Vc_k$ be given by
\[
\nu_t^*(u') \ = \
\begin{cases}
k+1, &\quad R_t^{\beta,n,k}(u')\leq0, \\
1, &\quad R_t^{\beta,n,k}(u')>0.
\end{cases}
\]
Then, from \eqref{BSDEnk_xi_nu3}, using also \eqref{Ineq_xi} and \eqref{Ineq_nu}, we obtain
\begin{equation}
\label{SaddlePoint}
G_T(\xi,\nu^*) \ \leq \ Y_0^{\beta,n,k} \ = \ G_T(\xi^\eps,\nu^*) + \eps\delta_T(\xi^\eps,\nu^*) \ \leq \ G_T(\xi^\eps,\nu) + \eps\delta_T(\xi^\eps,\nu),
\end{equation}
for any $\xi\in\Xi_n$ and $\nu\in\Vc_k$, where
\begin{align*}
\delta_T(\xi,\nu) \ &= \ \sum_{j=1}^m\int_0^T e^{-\beta t} \E^{\xi,\nu}\Big[\big(c(X_t^{x,i,u},I_{t^-}^i,j) - L_t^{\beta,n,k}(j)\big)1_{\{-1 \leq L_t^{\beta,n,k}(j) - c(X_t^{x,i,u},I_{t^-}^i,j) < 0\}} \\
&\quad \ + 1_{\{L_t^{\beta,n,k}(j) - c(X_t^{x,i,u},I_{t^-}^i,j) < -1\}}\Big] dt.
\end{align*}
Notice that $0\leq\delta_T(\xi,\nu)\leq \sum_{j=1}^m\int_0^T e^{-\beta t}dt \leq m/\beta$, therefore from \eqref{SaddlePoint} we have
\[
G_T(\xi,\nu^*) \ \leq \ Y_0^{\beta,n,k} \ \leq \ G_T(\xi^\eps,\nu) + \eps\frac{m}{\beta}, \qquad \forall\,\xi\in\Xi_n,\,\nu\in\Vc_k,
\]
which implies
\[
\inf_{\nu\in\Vc_k}\sup_{\xi\in\Xi_n}G_T(\xi,\nu) \ \leq \ Y_0^{\beta,n,k} \ \leq \ \sup_{\xi\in\Xi_n}\inf_{\nu\in\Vc_k}G_T(\xi,\nu) + \eps\frac{m}{\beta}.
\]
Since $\eps$ is arbitrary and the inequality $\sup_{\xi\in\Xi_n}\inf_{\nu\in\Vc_k}G_T(\xi,\nu)\leq\inf_{\nu\in\Vc_k}\sup_{\xi\in\Xi_n}G_T(\xi,\nu)$ holds, we deduce identities \eqref{DualFormula_T}. We also observe that, from \eqref{SaddlePoint} and the boundedness of $\delta_T(\xi^\eps,\nu^*)$, we have $Y_0^{\beta,n,k}=\lim_{\eps\rightarrow0^+}G_T(\xi^\eps,\nu^*)$, therefore $(\xi^\eps,\nu^*)$ in an $\eps$-saddle point of $G_T(\xi,\nu)$ on $\Xi_n\times\Vc_k$.

To obtain the dual representation formula \eqref{DualFormulaV^beta,n,k}, we now pass to the limit in \eqref{DualFormula_T} as $T\rightarrow\infty$. Firstly, we notice that, by estimates \eqref{EstimateYbeta,n,k} and \eqref{EstimateX_diss_xi,nu}, it follows from Lebesgue's dominated convergence theorem that
\begin{equation}
\label{DualConvergence1}
\sup_{\xi\in\Xi_n,\,\nu\in\Vc_k}\E^{\xi,\nu}\big[e^{-\beta T} |Y_T^{\beta,n,k}|\big] \ \overset{T\rightarrow\infty}{\longrightarrow} \ 0.
\end{equation}
Similarly, from the uniform linear growth condition of $f$ and $c$ with respect to $x$, the boundedness of $\xi\in\Xi_n$, and estimate \eqref{EstimateX_diss_xi,nu}, we have
\begin{align}
\sup_{\xi\in\Xi_n,\,\nu\in\Vc_k}\E^{\xi,\nu}\bigg[\int_T^\infty e^{-\beta t} |f(X_t^{x,i,u},I_t^i,\Gamma_t^u)| dt\bigg] \ &\overset{T\rightarrow\infty}{\longrightarrow} \ 0, \label{DualConvergence2} \\
\sup_{\xi\in\Xi_n,\,\nu\in\Vc_k}\E^{\xi,\nu}\bigg[\sum_{j=1}^m\int_T^\infty e^{-\beta t} c(X_t^{x,i,u},I_{t^-}^i,j) \xi_t(j) dt\bigg] \ &\overset{T\rightarrow\infty}{\longrightarrow} \ 0. \label{DualConvergence3}
\end{align}
The above convergence results imply that we can pass to the limit in \eqref{DualFormula_T} as $T\rightarrow\infty$, and we find (we are interested only in the first identity in \eqref{DualFormula_T})
\[
Y_0^{\beta,n,k} = \sup_{\xi\in\Xi_n}\inf_{\nu\in\Vc_k}\E^{\xi,\nu}\bigg[\int_0^\infty e^{-\beta t} f(X_t^{x,i,u},I_t^i,\Gamma_t^u) dt - \sum_{j=1}^m\int_0^\infty e^{-\beta t} c(X_t^{x,i,u},I_{t^-}^i,j) \xi_t(j) dt\bigg].
\]
The second integral in the above identity can be written with respect to $\pi$, using the definition of $\P^{\xi,\nu}$-compensator of $\pi$ (see, e.g., Theorem II.1.8(i) in \cite{jacod_shiryaev}), so that we end up with the dual representation formula \eqref{DualFormulaV^beta,n,k} for $V^{\beta,n,k}(x,i,u)=Y_0^{\beta,n,k}$.

\vspace{1mm}

\noindent\textbf{Step II.} \emph{Dual representation formula for $V^{\beta,n}$.} From the convergence $Y_0^{\beta,n,k}\searrow Y_0^{\beta,n,x,i,u}$ $\P$-a.s., as $k\rightarrow\infty$, we obtain
\begin{equation}
\label{DualInequality}
Y_0^{\beta,n,x,i,u} \ = \ \lim_{k\rightarrow\infty} \Big(\sup_{\xi\in\Xi_n}\inf_{\nu\in\Vc_k} G_\infty(\xi,\nu)\Big) \ \geq \ \sup_{\xi\in\Xi_n}\inf_{\nu\in\Vc} G_\infty(\xi,\nu),
\end{equation}
where
\[
G_\infty(\xi,\nu) \ = \ \E^{\xi,\nu}\bigg[\int_0^\infty e^{-\beta t}f(X_t^{x,i,u},I_t^i,\Gamma_t^u) dt - \int_0^\infty\int_{\I_m} e^{-\beta t} c(X_t^{x,i,u},I_{t^-}^i,j) \pi(dt,d\,\!j) \bigg].
\]
On the other hand, from \eqref{SaddlePoint} we have (recalling that $\delta_T(\xi,\nu)\leq m/\beta$)
\[
Y_0^{\beta,n,x,i,u} \ \leq \ Y_0^{\beta,n,k} \ \leq \ G_T(\xi^\eps,\nu) + \eps\frac{m}{\beta}, \qquad \forall\,\nu\in\Vc_k.
\]
Using the convergence results \eqref{DualConvergence1}-\eqref{DualConvergence2}-\eqref{DualConvergence3}, we can pass to the limit as $T\rightarrow\infty$ in the above inequalities, so that we obtain
\[
Y_0^{\beta,n,x,i,u} \ \leq \ Y_0^{\beta,n,k} \ \leq \ G_\infty(\xi^\eps,\nu) + \eps\frac{m}{\beta}, \qquad \forall\,\nu\in\Vc_k,
\]
which implies
\[
Y_0^{\beta,n,x,i,u} \ \leq \ \inf_{\nu\in\Vc_k} G_\infty(\xi^\eps,\nu)  + \eps\frac{m}{\beta}.
\]
Since $k$ and $\eps$ are arbitrary, we end up with
\[
Y_0^{\beta,n,x,i,u} \ \leq \ \sup_{\xi\in\Xi_n}\inf_{\nu\in\Vc} G_\infty(\xi,\nu).
\]
The above inequality, together with \eqref{DualInequality}, yields the thesis, recalling from Proposition \ref{T:Feynman-Kac}(ii) that $V^{\beta,n}(x,i)=Y_0^{\beta,n,x,i,u}$, for any $(x,i,u)\in\R^d\times\I_m\times\mathring U$.
\ep

\subsubsection{Estimates}

From the dual representation formula \eqref{DualFormulaV^beta,n}, we deduce the following estimates for $V^{\beta,n}$.

\begin{Corollary}
\label{C:EstimatesV^beta,n}
Let Assumptions {\bf (H1)}, {\bf (H2)}, {\bf (H3)}, and {\bf (H$U$)} hold. \begin{itemize}
\item[\textup{(i)}] For any $\beta>0$ and $n\in\N\backslash\{0\}$, we have
\begin{equation}
\label{V^beta,n_Lipschitz}
|V^{\beta,n}(x,i) - V^{\beta,n}(x',i)| \ \leq \ \frac{L_2}{\gamma}|x - x'|,
\end{equation}
for all $x,x'\in\R^d$, $i\in\I_m$.
\item[\textup{(ii)}] There exists a constant $\hat C\geq0$, depending only on $L_2$, $M:=\sup_{(i,u)\in\I_m\times U}|f(0,i,u)|$, and the constant $\bar C$ appearing in estimate \eqref{EstimateX_diss_xi,nu}, such that, for any $\beta>0$ and $n\in\N\backslash\{0\}$,
\begin{equation}
\label{V^beta,n_Bound}
|\beta V^{\beta,n}(x,i)| \ \leq \ \hat C(1 + |x|),
\end{equation}
for all $(x,i)\in\R^d\times\I_m$.
\item[\textup{(iii)}] For any $\beta>0$, $n\in\N\backslash\{0\}$, and $(x,i)\in\R^d\times\I_m$, we have
\begin{equation}
\label{V^beta,n_increasing}
V^{\beta,n}(x,i) \ \leq \ V^{\beta,n+1}(x,i).
\end{equation}
\end{itemize}
\end{Corollary}
\textbf{Proof.}
(i) \emph{Estimate \eqref{V^beta,n_Lipschitz}.} Fix $x,x'\in\R^d$, $i\in\I_m$, and $u\in\mathring U$. From the dual representation formula \eqref{DualFormulaV^beta,n} and since $c=c(i,j)$ does not depend on $x$ under {\bf (H3)}, we have
\[
|V^{\beta,n}(x,i) - V^{\beta,n}(x',i)| \ \leq \ \sup_{\xi\in\Xi_n,\,\nu\in\Vc} \E^{\xi,\nu}\bigg[\int_0^\infty e^{-\beta t}\big|f(X_t^{x,i,u},I_t^i,\Gamma_t^u) - f(X_t^{x',i,u},I_t^i,\Gamma_t^u)\big| dt\bigg].
\]
Then, using the Lipschitz property of $f$ in {\bf (H1)}(ii), together with estimate \eqref{EstimateX-X'}, we obtain
\[
|V^{\beta,n}(x,i) - V^{\beta,n}(x',i)| \ \leq \ L_2|x - x'|\int_0^\infty e^{-(\beta + \gamma) t} dt \ = \ \frac{L_2}{\beta + \gamma}|x - x'| \ \leq \ \frac{L_2}{\gamma}|x - x'|,
\]
from which the claim follows.

\vspace{1mm}

\noindent(ii) \emph{Estimate \eqref{V^beta,n_Bound}.} Given $(x,i)\in\R^d\times\I_m$ and $u\in\mathring U$, using again the dual representation formula \eqref{DualFormulaV^beta,n} we find (from the nonnegativity of $c$)
\begin{equation}
\label{ProofbetaVbeta}
\beta V^{\beta,n}(x,i) \ \leq \ \sup_{\xi\in\Xi_n,\,\nu\in\Vc} \E^{\xi,\nu}\bigg[\int_0^\infty \beta e^{-\beta t}\big|f(X_t^{x,i,u},I_t^i,\Gamma_t^u)\big| dt\bigg].
\end{equation}
On the other hand, take $\eps>0$ and define
\[
\xi^\eps \ \equiv \ \frac{\eps\beta}{m}\frac{1}{1 + \max_{i,j\in\I_m}c(i,j)}.
\]
We see that, if $\eps$ is small enough so that $\xi^\eps\leq n$, then $\xi^\eps\in\Xi_n$. Therefore, from \eqref{DualFormulaV^beta,n} we obtain (recall that the stochastic integral in \eqref{DualFormulaV^beta,n} with respect to $\pi$ can be written with respect to its $\P^{\xi^\eps,\nu}$-compensator)
\begin{align*}
V^{\beta,n}(x,i) \ &\geq \ \inf_{\nu\in\Vc}\E^{\xi^\eps,\nu}\bigg[\int_0^\infty e^{-\beta t}f(X_t^{x,i,u},I_t^i,\Gamma_t^u) dt - \sum_{j=1}^m\int_{\I_m} e^{-\beta t} c(X_t^{x,i,u},I_{t^-}^i,j) \xi^\eps(j) dt \bigg] \\
&\geq \ - \sup_{\nu\in\Vc}\E^{\xi^\eps,\nu}\bigg[\int_0^\infty e^{-\beta t}|f(X_t^{x,i,u},I_t^i,\Gamma_t^u)| dt\bigg] - \eps \\
&\geq \ - \sup_{\xi\in\Xi_n,\,\nu\in\Vc}\E^{\xi,\nu}\bigg[\int_0^\infty e^{-\beta t}|f(X_t^{x,i,u},I_t^i,\Gamma_t^u)| dt\bigg] - \eps.
\end{align*}
By the arbitrariness of $\eps$ and inequality \eqref{ProofbetaVbeta}, we conclude that
\[
|\beta V^{\beta,n}(x,i)| \ \leq \ \sup_{\xi\in\Xi_n,\,\nu\in\Vc} \E^{\xi,\nu}\bigg[\int_0^\infty \beta e^{-\beta t}\big|f(X_t^{x,i,u},I_t^i,\Gamma_t^u)\big| dt\bigg].
\]
From the inequality $f(x,i,u)\leq L_2|x| + M$ and estimate \eqref{EstimateX_diss_xi,nu}, we see that there exists some positive constant $\hat C$ (only depending on $L_2$, $M$, $\bar C$) such that
\[
|\beta V^{\beta,n}(x,i)| \ \leq \ \hat C(1 + |x|) \int_0^\infty \beta e^{-\beta t} dt \ = \ \hat C(1 + |x|).
\]

\vspace{1mm}

\noindent(iii) \emph{Monotone property  \eqref{V^beta,n_increasing}.} Inequality \eqref{V^beta,n_increasing} follows from the dual representation formula \eqref{DualFormulaV^beta,n} for $V^{\beta,n}$, noting that $\Xi_n\subset\Xi_{n+1}$, $\forall\,n\in\N\backslash\{0\}$.
\ep

\subsubsection{Convergence of $V^{\beta,n}$ towards $V^\beta$}

From Corollary \ref{C:EstimatesV^beta,n} it follows that, for any $\beta>0$, the sequence $(V^{\beta,n})_{n\geq1}$ is monotone nondecreasing and satisfies the following linear growth condition:
\begin{equation}
\label{V^beta,n_linear_growth}
V^{\beta,n}(x,i) \ \leq \ \frac{\hat C}{\beta} + \frac{L_2}{\gamma}|x|,
\end{equation}
for all $n\in\N\backslash\{0\}$ and $(x,i)\in\R^d\times\I_m$. As a consequence, we have the following result.

\begin{Proposition}
\label{P:V^beta}
Let Assumptions {\bf (H1)}, {\bf (H2)}, {\bf (H3)}, and {\bf (H$U$)} hold. Then, for any $\beta>0$, there exists a function $v^\beta\colon\R^d\times\I_m\rightarrow\R$ given by $v^\beta(x,i) = \lim_{n\rightarrow\infty}V^{\beta,n}(x,i)$, for all $(x,i)\in\R^d\times\I_m$, which admits the dual representation formula
\begin{align*}
v^\beta(x,i) \ = \ \sup_{\xi\in\Xi}\inf_{\nu\in\Vc}\E^{\xi,\nu}\bigg[&\int_0^\infty e^{-\beta t}f(X_t^{x,i,u},I_t^i,\Gamma_t^u) dt \\
& - \int_0^\infty\int_{\I_m} e^{-\beta t} c(I_{t^-}^i,j) \pi(dt,d\,\!j) \bigg], \qquad \forall\,(x,i)\in\R^d\times\I_m,
\end{align*}
for any $u\in\mathring U$. Moreover, $v^\beta$ has the following properties:
\begin{align*}
|v^\beta(x,i) - v^\beta(x',i)| \ &\leq \ \frac{L_2}{\gamma}|x - x'|, \\
|\beta v^\beta(x,i)| \ &\leq \ \hat C(1 + |x|),
\end{align*}
for all $x,x'\in\R^d$ and $i\in\I_m$, where the constant $\hat C$ is the same as in Corollary \ref{C:EstimatesV^beta,n}(ii).
\end{Proposition}
\textbf{Proof.}
The existence of $v^\beta$ follows from the monotone property of $(V^{\beta,n})_{n\geq1}$ (which is a consequence of \eqref{V^beta,n_increasing}) and from the uniform linear growth condition \eqref{V^beta,n_linear_growth}. The dual representation formula holds since $V^{\beta,n}$ satisfies \eqref{DualFormulaV^beta,n} and $\Xi_n\subset\Xi_{n+1}\subset\cdots\subset\cup_{n\geq1}\Xi_n=\Xi$. Finally, the two stated properties of $v^\beta$ are implied by estimates \eqref{V^beta,n_Lipschitz} and \eqref{V^beta,n_Bound}.
\ep

\begin{Proposition}
\label{P:VbetaViscSol}
Let Assumptions {\bf (H1)}, {\bf (H2)}, {\bf (H3)}, and {\bf (H$U$)} hold. Then, the function $v^\beta$ is the unique continuous viscosity solution to system \eqref{HJB_ellipticintro} satisfying a linear growth condition
\[
\sup_{(x,i)\in\R^d\times\I_m} \frac{|v^\beta(x,i)|}{1 + |x|} \ < \ \infty.
\]
In particular, $v^\beta$ coincides with the function $V^\beta$ defined by \eqref{V^beta}. Moreover, we have
\begin{equation}
\label{Vbetai-j}
|v^\beta(x,i) - v^\beta(x,j)| \ \leq \ \hat c \ := \ \max_{i,j\in\I_m}c(i,j),
\end{equation}
for all $x\in\R^d$ and $i,j\in\I_m$.
\end{Proposition}
\begin{Remark}
\label{R:DualProblem}
{\rm
We refer to the function $v^\beta$ introduced in Proposition \ref{P:V^beta} as \emph{dual value function} of the \emph{dual robust switching control problem}. On the other hand, we say that $V^\beta$ is the \emph{primal value function} and the associated infinite horizon robust feedback switching control, recalled in Section \ref{S:Motivation}, is the \emph{primal robust switching control problem}. Then, Propositions \ref{P:V^beta} and \ref{P:VbetaViscSol} state that these two value functions coincide. Notice that the dual control problem is a two-player zero-sum stochastic differential game of the type \emph{control vs control}. In general, the lower and upper value functions of a control vs control (more precisely, \emph{open-loop control vs open-loop control}) game do not coincide with those associated to the Elliott-Kalton version of the game, where the \emph{strategy vs control} formulation is adopted (see, e.g., Exercise 2.1(ii), Chapter VIII, in \cite{bardi_capuzzo-dolcetta}). However, we observe that our dual control problem is in weak form as in \cite{pham_zhang}, therefore, as emphasized in \cite{pham_zhang}, it might be interpreted as a \emph{feedback control vs feedback control} game.
}
\ep
\end{Remark}
\textbf{Proof (of Proposition \ref{P:VbetaViscSol}).}
Recall that $V^{\beta,n}$ satisfies \eqref{V^beta,n_Lipschitz} and $v^\beta$ is the pointwise limit of the sequence $(V^{\beta,n})_{n\geq1}$. Then, we have
\begin{equation}
\label{V^beta=*0}
|V^{\beta,n}(x',i) - v^\beta(x,i)| \ \leq \ \frac{L_2}{\gamma}|x - x'| + |V^{\beta,n}(x,i) - v^\beta(x,i)| \ \underset{\substack{n\rightarrow\infty\\ x'\rightarrow x}}{\longrightarrow} \ 0
\end{equation}
Therefore, for all $(x,i)\in\R^d\times\I_m$,
\begin{equation}
\label{V^beta=*}
v^\beta(x,i) \ = \ \lim_{\substack{n\rightarrow\infty\\ x'\rightarrow x}} V^{\beta,n}(x',i) \ = \ \liminf_{n\rightarrow\infty}{\!}_*\, V^{\beta,n}(x,i) \ = \ \limsup_{n\rightarrow\infty}{\!}^*\, V^{\beta,n}(x,i),
\end{equation}
where
\[
\liminf_{n\rightarrow\infty}{\!}_*\, V^{\beta,n}(x,i) \ := \ \liminf_{\substack{n\rightarrow\infty\\ x'\rightarrow x}} V^{\beta,n}(x',i), \qquad \limsup_{n\rightarrow\infty}{\!}^*\, V^{\beta,n}(x,i) \ := \ \limsup_{\substack{n\rightarrow\infty\\ x'\rightarrow x}} V^{\beta,n}(x',i).
\]

\vspace{1mm}

\noindent\textbf{Step I.} \emph{Viscosity supersolution property of $v^\beta$ to \eqref{HJB_ellipticintro}.} Let $(x,i)\in\R^d\times\I_m$ and $(p,M)\in J^{2,-}v^\beta(x,i)$ (the second-order subjet of $v^\beta$ at $(x,i)$, see, e.g., Section 2 in \cite{crandall_ishii_lions}). From \eqref{V^beta=*} and Lemma 6.1 (see also Remark 6.2) in \cite{crandall_ishii_lions}, it follows that we can find the following sequences
\[
n_k \ \overset{k\rightarrow\infty}{\longrightarrow} \ \infty, \qquad x_k\in\R^d, \qquad (p_k,M_k)\in J^{2,-} V^{\beta,n_k}(x_k,i),
\]
satisfying
\begin{equation}
\label{Convergences}
(x_k,V^{\beta,n_k}(x_k,i),p_k,M_k) \ \overset{k\rightarrow\infty}{\longrightarrow} \ (x,v^\beta(x,i),p,M).
\end{equation}
For any $j\in\I_m$, we also have the convergence $V^{\beta,n_k}(x_k,j)\rightarrow v^\beta(x,j)$, which can be proved proceeding as in \eqref{V^beta=*0} with $j$ in place of $i$. From the supersolution property of $V^{\beta,n_k}$ to \eqref{HJB_penalized} stated in Proposition \ref{T:Feynman-Kac}(ii), we have
\begin{align}
\label{HJB_penalized_Proof}
\beta V^{\beta,n_k}(x_k,i) - \inf_{u\in U}\Big[b(x_k,i,u).p_k + \frac{1}{2}\text{tr}\big(\sigma\sigma\trans(x_k,i,u)M_k\big) + f(x_k,i,u)\Big]\;\;\,& \\
- n_k\sum_{j=1}^m \big[V^{\beta,n_k}(x_k,j) - V^{\beta,n_k}(x_k,i)
- c(i,j)\big]^+ &\geq \, 0. \notag
\end{align}
Let us prove that
\begin{equation}
\label{v^beta_Mv^beta}
v^\beta(x,i) - \max_{j\neq i}\big[v^\beta(x,j) - c(i,j)\big] \ \geq \ 0.
\end{equation}
On the contrary, suppose that there exists some $j_0\in\I_m$, $j_0\neq i$, such that
\[
v^\beta(x,i) - v^\beta(x,j_0) + c(i,j_0) \ < \ 0.
\]
From the convergences $V^{\beta,n_k}(x_k,i)\rightarrow v^\beta(x,i)$, $V^{\beta,n_k}(x_k,j_0)\rightarrow v^\beta(x,j_0)$, it follows that there exist $\eps>0$ and $k_\eps\in\N$ such that
\[
V^{\beta,n_k}(x_k,i) - V^{\beta,n_k}(x_k,j_0) + c(i,j_0) \ \leq \ - \eps, \qquad \forall\,k\geq k_\eps.
\]
As a consequence, we have
\[
\sum_{j=1}^m \big[V^{\beta,n_k}(x_k,j) - V^{\beta,n_k}(x_k,i)
- c(i,j)\big]^+ \ \geq \ \eps, \qquad \forall\,k\geq k_\eps.
\]
Letting $k\rightarrow\infty$ into \eqref{HJB_penalized_Proof}, we find a contradiction, so that \eqref{v^beta_Mv^beta} holds. On the other hand, from \eqref{HJB_penalized_Proof} we have
\[
\beta V^{\beta,n_k}(x_k,i) - \inf_{u\in U}\Big[b(x_k,i,u).p_k + \frac{1}{2}\text{tr}\big(\sigma\sigma\trans(x_k,i,u)M_k\big) + f(x_k,i,u)\Big] \ \geq \ 0.
\]
Sending $k\rightarrow\infty$, using \eqref{Convergences} and the continuity of $b$, $\sigma$, $f$, we conclude that
\[
\beta v^\beta(x,i) - \inf_{u\in U}\Big[b(x,i,u).p + \frac{1}{2}\text{tr}\big(\sigma\sigma\trans(x,i,u)M\big) + f(x,i,u)\Big] \ \geq \ 0.
\]

\vspace{1mm}

\noindent\textbf{Step II.} \emph{Viscosity subsolution property of $v^\beta$ to \eqref{HJB_ellipticintro}.} Let $(x,i)\in\R^d\times\I_m$ and $(p,M)\in J^{2,+} v^\beta(x,i)$ (the second-order superjet of $v^\beta$ at $(x,i)$, see Section 2 in \cite{crandall_ishii_lions}) such that
\begin{equation}
\label{v^beta_Mv^beta2}
v^\beta(x,i) - \max_{j\neq i}\big[v^\beta(x,j) - c(i,j)\big] \ > \ 0.
\end{equation}
From \eqref{V^beta=*} and Lemma 6.1 in \cite{crandall_ishii_lions}, we see that we can find the following sequences
\[
n_k \ \overset{k\rightarrow\infty}{\longrightarrow} \ \infty, \qquad x_k\in\R^d, \qquad (p_k,M_k)\in J^{2,+} V^{\beta,n_k}(x_k,i),
\]
satisfying
\[
(x_k,V^{\beta,n_k}(x_k,i),p_k,M_k) \ \overset{k\rightarrow\infty}{\longrightarrow} \ (x,v^\beta(x,i),p,M).
\]
Moreover, for any $j\in\I_m$, we also have the convergence $V^{\beta,n_k}(x_k,j)\rightarrow v^\beta(x,j)$. Using the subsolution property of $V^{\beta,n_k}$ to \eqref{HJB_penalized}, we find
\begin{align}
\label{HJB_penalized_Proof2}
\beta V^{\beta,n_k}(x_k,i) - \inf_{u\in U}\Big[b(x_k,i,u).p_k + \frac{1}{2}\text{tr}\big(\sigma\sigma\trans(x_k,i,u)M_k\big) + f(x_k,i,u)\Big]\;\;\,& \\
- n_k\sum_{j=1}^m \big[V^{\beta,n_k}(x_k,j) - V^{\beta,n_k}(x_k,i)
- c(i,j)\big]^+ &\leq \, 0. \notag
\end{align}
From \eqref{v^beta_Mv^beta2} and the convergence $V^{\beta,n_k}(x_k,j)\rightarrow v^\beta(x,j)$, $\forall\,j\in\I_m$, we see that there exist $\eps>0$ and $k_\eps\in\N$ such that
\[
V^{\beta,n_k}(x_k,i) - \max_{j\neq i}\big[V^{\beta,n_k}(x_k,j) - c(i,j)\big] \ \geq \ \eps, \qquad \forall\,k\geq k_\eps.
\]
Therefore, for all $k\geq k_\eps$, we have
\[
\sum_{j=1}^m \big[V^{\beta,n_k}(x_k,j) - V^{\beta,n_k}(x_k,i)
- c(i,j)\big]^+ \ = \ 0.
\]
Hence, letting $k\rightarrow\infty$ into \eqref{HJB_penalized_Proof2}, we end up with
\[
\beta v^\beta(x,i) - \inf_{u\in U}\Big[b(x,i,u).p + \frac{1}{2}\text{tr}\big(\sigma\sigma\trans(x,i,u)M\big) + f(x,i,u)\Big] \ \leq \ 0.
\]

\vspace{1mm}

\noindent\textbf{Step III.} \emph{Identification $v^\beta\equiv V^\beta$.} From Proposition \ref{P:VbetaStochasticPerron} we know that $V^\beta$ is the unique continuous viscosity solution to system \eqref{HJB_ellipticintro} satisfying a linear growth condition, so that the claim follows.

\vspace{1mm}

\noindent\textbf{Step IV.} \emph{Estimate \eqref{Vbetai-j}.} Finally, to prove estimate \eqref{Vbetai-j}, we notice that from the identification $v^\beta\equiv V^\beta$ and inequality \eqref{VbetaInequality}, we have
\[
V^{\beta}(x,i) \ \geq \ \max_{j\neq i}\big[V^\beta(x,j) - c(i,j)\big] \ \geq \ V^\beta(x,j) - c(i,j), \quad \forall\,(x,i,j)\in\R^d\times\I_m^2,\,j\neq i.
\]
This implies that
\[
V^\beta(x,j) - V^\beta(x,i) \ \leq \ c(i,j) \ \leq \ \max_{i,j\in\I_m}c(i,j),
\]
from which estimate \eqref{Vbetai-j} follows.
\ep

\subsection{Convergence results for $V^\beta$}

We are now in a position to study the asymptotic behavior of $V^\beta$. More precisely, we have the following result, which proves all the statements of Theorem \ref{T:Main} concerning $V^\beta$ and the existence of a viscosity solution to the ergodic system \eqref{HJB_ergodic}.

\begin{Proposition} \label{P:convergence}
Let Assumptions {\bf (H1)}, {\bf (H2)}, {\bf (H3)}, and {\bf (H$U$)} hold. Then,  there exists a viscosity solution $(\lambda,\phi)$, with $\phi(\cdot,i)$ Lipschitz, for any $i\in\I_m$, and $\phi(0,i_0)=0$ for some fixed $i_0\in\I_m$, to the ergodic system \eqref{HJB_ergodic}, such that
\begin{align*}
\beta V^\beta(x,i)& \ \,\,\,\overset{\beta\rightarrow0^+}{\longrightarrow} \ \lambda, \qquad\qquad\, \forall\,(x,i)\in\R^d\times\I_m, \\
V^{\beta_k}(\cdot,i) - V^{\beta_k}(0,i_0)& \ \underset{\text{in }C(\R^d)}{\overset{k\rightarrow\infty}{\longrightarrow}} \ \phi(\cdot,i), \qquad \forall\,i\in\I_m,
\end{align*}
for some sequence $(\beta_k)_{k\in\N}$, with $\beta_k\searrow0^+$.
\end{Proposition}
\textbf{Proof.}
Fix $i_0\in\I_m$ and, for any $\beta>0$, set
\[
\lambda_i^\beta \ := \ \beta V^\beta(0,i), \qquad \phi^\beta(x,i) \ := \ V^\beta(x,i) - V^\beta(0,i_0),
\]
for all $(x,i)\in\R^d\times\I_m$. From Proposition \ref{P:V^beta}, estimate \eqref{Vbetai-j}, and the identification $v^\beta\equiv V^\beta$ stated in Proposition \ref{P:VbetaViscSol}, we have
\[
\sup_{\beta>0}|\lambda_i^\beta| \ \leq \ \hat C, \qquad \sup_{\beta>0}|\phi^\beta(x,i)| \ \leq \ \frac{L_2}{\gamma}|x| + \hat c.
\]
As a consequence, by classical arguments based on the Bolzano-Weierstrass Theorem and the Ascoli-Arzel\`a Theorem, see e.g. \cite{fuh-hu-tess}, we can find a sequence $(\beta_k)_{k\in\N}$, with 
$\beta_k\searrow0^+$ s.t. 
\[
\lambda_i^{\beta_k} \ \overset{k\rightarrow\infty}{\longrightarrow} \ \lambda_i, \qquad \phi^{\beta_k}(\cdot,i) \ \underset{\text{in }C(\R^d)}{\overset{k\rightarrow\infty}{\longrightarrow}} \ \phi(\cdot,i),
\]
for some $\lambda_i\in\R$ and $\phi\colon\R^d\times\I_m\rightarrow\R$ satisfying $|\phi(x,i)|\leq L_2|x|/\gamma+\hat c$, $|\phi(x,i)-\phi(x',i)|\leq L_2|x-x'|/\gamma$, and $\phi(0,i_0)=0$. Notice that, from estimate \eqref{Vbetai-j} we obtain
\[
\big|\lambda_i^{\beta_k} - \lambda_j^{\beta_k}\big| \ = \ \beta_k\big|V^{\beta_k}(0,i) - V^{\beta_k}(0,j)\big| \ \leq \ \beta_k\hat c \ \overset{k\rightarrow\infty}{\longrightarrow} \ 0,
\]
therefore $\lambda:=\lambda_i=\lambda_j$, for all $i,j\in\I_m$. More generally, we have
\[
\big|\beta_k V^{\beta_k}(x,i) - \lambda_j^{\beta_k}\big| \ = \ \beta_k\big|V^{\beta_k}(x,i) - V^{\beta_k}(0,j)\big| \ \leq \ \beta_k\frac{L_2}{\gamma}|x| + \beta_k\hat c \ \overset{k\rightarrow\infty}{\longrightarrow} \ 0,
\]
which implies that
\[
\beta_k V^{\beta_k}(x,i) \ \overset{k\rightarrow\infty}{\longrightarrow} \ \lambda, \qquad \forall\,(x,i)\in\R^d\times\I_m.
\]
We now prove that $(\lambda,\phi)$ is a viscosity solution to the ergodic system \eqref{HJB_ergodic}. To this end, we begin noting that, from the viscosity properties of $V^\beta$ stated in Proposition \ref{P:VbetaViscSol} (see also Proposition \ref{P:VbetaStochasticPerron}), it follows that, for any $i\in\I_m$, $\phi^\beta(\cdot,i)$ is a viscosity solution to the following elliptic equation:
\begin{align*}
\min\Big\{\lambda_{i_0}^\beta + \beta \phi^\beta(x,i) - \inf_{u\in U}\big[\Lc^{i,u}\phi^\beta(x,i) + f(x,i,u)\big],& \\
\phi^\beta(x,i) - \max_{j\neq i}\big[\phi^\beta(x,j) - c(i,j)\big]\Big\}& \ = \ 0, \qquad \forall\,x\in\R^d.
\end{align*}
Then, we define
\begin{align*}
F_k(x,i,r,p,M) \ := \ \min\Big\{&\lambda_{i_0}^{\beta_k} + \beta_k r - \inf_{u\in U}\big[b(x,i,u).p + \frac{1}{2}\text{tr}\big(\sigma\sigma\trans(x,i,u)M\big) + f(x,i,u)\big], \\
& r - \max_{j\neq i}\big[\phi^{\beta_k}(x,j) - c(i,j)\big]\Big\},
\end{align*}
for all $k\in\N$, and
\begin{align*}
F_\infty(x,i,r,p,M) \ := \ \min\Big\{&\lambda - \inf_{u\in U}\big[b(x,i,u).p + \frac{1}{2}\text{tr}\big(\sigma\sigma\trans(x,i,u)M\big) + f(x,i,u)\big], \\
& r - \max_{j\neq i}\big[\phi(x,j) - c(i,j)\big]\Big\},
\end{align*}
for any $(x,i,r,p,M)\in\R^d\times\I_m\times\R\times\R^d\times\R^{d\times d}$. We see that
\[
\lim_{k\rightarrow\infty} F_k(x,i,r,p,M) \ = \ F_\infty(x,i,r,p,M).
\]
Then, from stability results of viscosity solutions (see, e.g., Lemma 6.1 and Remark 6.3 in \cite{crandall_ishii_lions}), we deduce that, for any $i\in\I_m$, the function $\phi(\cdot,i)$ is a viscosity solution to the elliptic equation
\[
F_\infty(x,i,\phi(x,i),D_x\phi(x,i),D_x^2\phi(x,i)) \ = \ 0, \qquad \forall\,x\in\R^d.
\]
As a consequence, we conclude that $(\lambda,\phi)$ is a viscosity solution to the ergodic system \eqref{HJB_ergodic}. Finally, we notice that the all family $(\beta V^\beta(x,i))_{\beta>0}$ converges to $\lambda$ as $\beta\rightarrow0^+$, since, as stated in Proposition \ref{P:V-->lambda}, $\lambda$ is uniquely determined.
\ep

\appendix

\setcounter{equation}{0} \setcounter{Assumption}{0}
\setcounter{Theorem}{0} \setcounter{Proposition}{0}
\setcounter{Corollary}{0} \setcounter{Lemma}{0}
\setcounter{Definition}{0} \setcounter{Remark}{0}

\renewcommand\thesection{Appendix}

\section{}

\renewcommand\thesection{\Alph{subsection}}

\renewcommand\thesubsection{\Alph{subsection}}

\subsection{Feynman-Kac formula}

The present appendix is devoted to the proof of Proposition \ref{T:Feynman-Kac}. Unfortunately, we did not find a reference for it in the literature. Indeed, even though Proposition 3.3 and Theorem 3.1 in \cite{CFP14} do almost the job, they do not apply to system \eqref{HJB_penalized} due to the presence of the nonlocal term. For this reason, in the present appendix we state the results, recalling only the main steps of their proofs, since they are very similar to those of Proposition 3.3 and Theorem 3.1 in \cite{CFP14}.

\subsubsection{Maximal solution to BSDE \eqref{BSDE}-\eqref{BSDE_Constraint}}

We begin addressing the problem of existence and uniqueness of the maximal solution (see Definition \ref{D:Maximal}) to the BSDE with partially nonnegative jumps on infinite horizon \eqref{BSDE}-\eqref{BSDE_Constraint}. Concerning uniqueness, we have the following result.

\begin{Lemma}
\label{L:Uniqueness}
Suppose that Assumption {\bf (H1)} holds. Then, for any $\beta>0$, $n\in\N$, $(x,i,u)\in\R^d\times\I_m\times U$, there exists at most one maximal solution to equation \eqref{BSDE}-\eqref{BSDE_Constraint}.
\end{Lemma}
\textbf{Proof.}
The uniqueness of the $Y$ component follows by definition. Now, consider two maximal solutions $(Y,Z,L,R,K)$, $(Y,Z',L',R',K')$ in ${\bf S_{\textup{loc}}^2}\times{\bf L_{\textup{loc}}^2(W)}\times{\bf L_{\textup{loc}}^2(\tilde\pi)}\times{\bf L_{\textup{loc}}^2(\tilde\mu)}\times{\bf K_{\textup{loc}}^2}$ to \eqref{BSDE}-\eqref{BSDE_Constraint}. Taking their difference, and identifying the Brownian and finite variation parts, we see that $Z=Z'$. Afterwards, recalling that the marked point processes associated to $\pi$ and $\mu$ have disjoint (due to the independence of $\pi$ and $\mu$) totally inaccessible jumps, while $K$ and $K'$ have predictable jumps, we conclude that $L=L'$ and $R=R'$, so that $K=K'$. For more details, we refer to Remark 3.1 in \cite{CFP14}.
\ep

\vspace{3mm}

The existence of the maximal solution to \eqref{BSDE}-\eqref{BSDE_Constraint} is based on a penalization method. More precisely, for any $\beta>0$, $n,k\in\N$, and $(x,i,u)\in\R^d\times\I_m\times U$, consider the following doubly indexed penalized backward stochastic differential equation on infinite horizon:
\beq
\label{BSDEnk}
Y_t^{\beta,n,k} &=& Y_T^{\beta,n,k} - \beta\int_t^T Y_s^{\beta,n,k}ds + \int_t^T f(X_s^{x,i,u},I_s^i,\Gamma_s^u) ds - \sum_{j=1}^m \int_t^T L_s^{\beta,n,k}(j) ds \\
&& + \; n\sum_{j=1}^m\int_t^T \!\! \big[L_s^{\beta,n,k}(j) - c(X_s^{x,i,u},I_{s^-}^i,j)\big]^+ ds - k\int_t^T \! \int_U \! \big[R_s^{\beta,n,k}(u')\big]^- \vartheta_\mu(du')ds \notag \\
&& - \; \int_t^T Z_s^{\beta,n,k} dW_s - \int_t^T\int_{\I_m} L_s^{\beta,n,k}(j) \tilde\pi(ds,d\,\!j) - \int_t^T\int_U R_s^{\beta,n,k}(u') \tilde\mu(ds,du'), \notag
\enq
for any $0\leq t\leq T$, $T\in[0,\infty)$, where $h^-=-\min(h,0)$ denotes the negative part of the function $h$. Then, we have the following result.

\begin{Lemma}
\label{L:PenalizedBSDE}
Let Assumptions {\bf (H1)} and {\bf (H2)} hold. Then, for any $\beta>0$, $n,k\in\N$, $(x,i,u)\in\R^d\times\I_m\times U$, there exists a solution $(Y^{\beta,n,k,x,i,u},Z^{\beta,n,k,x,i,u},L^{\beta,n,k,x,i,u},R^{\beta,n,k,x,i,u})$ $\in$ ${\bf S_{\textup{loc}}^2}\times{\bf L_{\textup{loc}}^2(W)}\times{\bf L_{\textup{loc}}^2(\tilde\pi)}\times{\bf L_{\textup{loc}}^2(\tilde\mu)}$ to \eqref{BSDEnk} such that
\begin{equation}
\label{EstimateYbeta,n,k}
\big|Y_t^{\beta,n,k,x,i,u}\big| \ \leq \ \frac{C_{b,\sigma,f}}{\beta} \big(1 + \big|X_t^{x,i,u}\big|\big), \qquad \forall\,t\geq0,
\end{equation}
where $C_{b,\sigma,f}$ is a positive constant, depending only on $b,\sigma,f$. This latter solution is unique among all quadruplets $(Y,Z,L,R)\in{\bf S_{\textup{loc}}^2}\times{\bf L_{\textup{loc}}^2(W)}\times{\bf L_{\textup{loc}}^2(\tilde\pi)}\times{\bf L_{\textup{loc}}^2(\tilde\mu)}$ such that, for some constant $C\geq0$ (possibly depending on $\beta,n,k,x,i,u$), we have $|Y_t|\leq C(1+|X_t^{x,i,u}|)$, for all $t\geq0$. 
\end{Lemma}
\textbf{Proof.}
The result follows from the same arguments as in the proof of Proposition 3.1 in \cite{CFP14}. Here, we give simply a sketch of the proof.\\
\emph{Uniqueness.} Consider two solutions $(Y,Z,L,R)$, $(Y',Z',L',R')$ to \eqref{BSDEnk}, satisfying a linear growth condition as stated in Lemma \ref{L:PenalizedBSDE}, and apply It\^o formula to the difference $e^{-2\beta s}|Y_s-Y_s'|^2$ between $s=t\geq0$ and $s=T\geq t$. Then, from the resulting expression we see that there exists $\xi\in\Xi_k$ such that, taking the expectation $\E^{1,\xi}$ with respect to $\P^{1,\xi}$ (we denote by $\P^{1,\xi}$ the probability measure $\P^{\nu,\xi}$ with $\nu\equiv1$), we obtain
\[
\E^{1,\xi}\big[|Y_t-Y_t'|^2\big] \ \leq \ e^{-2\beta(T-t)}\E^{1,\xi}\big[|Y_T-Y_T'|^2\big].
\]
Then, using the growth condition of $Y_T,Y_T'$ together with estimate \eqref{EstimateX_diss_xi,nu}, we see that $Y-Y'=0$. Finally, since $Y=Y'$, the identities $Z=Z'$, $L=L'$, $R=R'$ can be proved proceeding as in Lemma \ref{L:Uniqueness}.

\vspace{1mm}

\noindent\emph{Existence.} The proof consists in approximating equation \eqref{BSDEnk} through a sequence of BSDEs with finite time horizon and zero terminal condition. More precisely, for any $T>0$ and $(x,i,u)\in\R^d\times\I_m\times U$, we consider the following backward stochastic differential equation on $[0,T]$:
\begin{align}
\label{BSDEnk_T}
Y_t^{T,\beta,n,k} &= - \beta\int_t^T Y_s^{T,\beta,n,k}ds + \int_t^T f(X_s^{x,i,u},I_s^i,\Gamma_s^u) ds - \sum_{j=1}^m \int_t^T L_s^{T,\beta,n,k}(j) ds \\
& + n\sum_{j=1}^m\int_t^T \big[L_s^{T,\beta,n,k}(j) - c(X_s^{x,i,u},I_{s^-}^i,j)\big]^+ ds - k\int_t^T\int_U \big[R_s^{T,\beta,n,k}(u')\big]^- \vartheta_\mu(du')ds \notag \\
& - \int_t^T Z_s^{T,\beta,n,k} dW_s - \int_t^T\int_{\I_m} L_s^{T,\beta,n,k}(j) \tilde\pi(ds,d\,\!j) - \int_t^T\int_U R_s^{T,\beta,n,k}(u') \tilde\mu(ds,du'), \notag
\end{align}
for any $0\leq t\leq T$. From Lemma 2.4 in \cite{tang_li} we know that there exists a unique solution $(Y^{T,\beta,n,k,x,i,u},Z^{T,\beta,n,k,x,i,u},L^{T,\beta,n,k,x,i,u},R^{T,\beta,n,k,x,i,u})\in{\bf S_T^2}\times{\bf L_T^2(W)}\times{\bf L_T^2(\tilde\pi)}\times{\bf L_T^2(\tilde\mu)}$ to equation \eqref{BSDEnk_T}. Proceeding as in Proposition 3.1 in \cite{CFP14}, we then exploit Girsanov's theorem to investigate the differences $Y^T-Y^{T'}$, $Z^T-Z^{T'}$, $L^T-L^{T'}$, $U^T-U^{T'}$. In particular, we are able to determine opportune estimates for those differences, which allow to pass to the limit as $T\rightarrow\infty$ in equation \eqref{BSDEnk_T} and to end up with the solution to \eqref{BSDEnk}.
\ep

\vspace{3mm}

We can now state the following existence and uniqueness result for the BSDE with partially nonnegative jumps on infinite horizon \eqref{BSDE}-\eqref{BSDE_Constraint}, which in particular proves statement (i) of Proposition \ref{T:Feynman-Kac}.

\begin{Proposition}
\label{P:PenalizedBSDE}
Let Assumptions {\bf (H1)} and {\bf (H2)} hold. Then, for any $\beta>0$, $n\in\N$, $(x,i,u)\in\R^d\times\I_m\times U$, there exists a unique maximal solution $(Y^{\beta,n,x,i,u},Z^{\beta,n,x,i,u},L^{\beta,n,x,i,u},$ $R^{\beta,n,x,i,u},K^{\beta,n,x,i,u})\in{\bf S_{\textup{loc}}^2}\times{\bf L_{\textup{loc}}^2(W)}\times{\bf L_{\textup{loc}}^2(\tilde\pi)}\times{\bf L_{\textup{loc}}^2(\tilde\mu)}\times{\bf K_{\textup{loc}}^2}$ to \eqref{BSDE}-\eqref{BSDE_Constraint} such that:
\begin{itemize}
\item[\textup{(i)}] For all $t\geq0$, $Y_t^{\beta,n,k,x,i,u}\searrow Y_t^{\beta,n,x,i,u}$ $\P$-a.s., as $k\rightarrow\infty$.
\item[\textup{(ii)}] For all $T>0$, $(Z_{|[0,T]}^{\beta,n,k,x,i,u},L_{|[0,T]}^{\beta,n,k,x,i,u},R_{|[0,T]}^{\beta,n,k,x,i,u})_k$ strongly $($resp. weakly$)$ converges to $(Z_{|[0,T]}^{\beta,n,x,i,u},L_{|[0,T]}^{\beta,n,x,i,u},R_{|[0,T]}^{\beta,n,x,i,u})$ in ${\bf L_T^p(W)}\times{\bf L_T^p(\tilde\pi)}\times{\bf L_T^p(\tilde\mu)}$, for any $p\in[1,2)$ $($resp. in ${\bf L_T^2(W)}\times{\bf L_T^2(\tilde\pi)}\times{\bf L_T^2(\tilde\mu)}$$)$.
\item[\textup{(iii)}] For all $t\geq0$, $(K_t^{\beta,n,k,x,i,u})_k$ weakly converges to $K_t^{\beta,n,x,i,u}$ in ${\bf L^2(}\Omega,\Fc_t,\P{\bf)}$.
\end{itemize}
\end{Proposition}
\textbf{Proof.}
The proof can be done proceeding as in Proposition 3.3 in \cite{CFP14}. We just recall the main steps. Firstly, for all $t\geq0$, the nonincreasing property of the sequence $(Y_t^{\beta,n,k,x,i,u})_k$ follows from the comparison theorem for BSDEs with jumps, see, e.g., Theorem 4.2 in \cite{quenez_sulem}. The monotonicity property provides the existence of a limiting adapted process $Y^{\beta,n,x,i,u}$ satisfying estimate \eqref{EstimateYbeta,n,k}. Afterwards, for any $T>0$, we consider equation \eqref{BSDE}-\eqref{BSDE_Constraint} on $[0,T]$ with terminal condition $Y_T^{\beta,n,x,i,u}$. Then, from Theorem 2.1 in \cite{KP12} it follows that there exists a unique maximal solution to \eqref{BSDE}-\eqref{BSDE_Constraint} (Theorem 2.1 in \cite{KP12} applies to minimal solutions; however, simply notice that if $Y$ is a maximal solution to \eqref{BSDE}-\eqref{BSDE_Constraint}, then $-Y$ is a minimal solution to a certain BSDE with partially nonpositive jumps to which Theorem 2.1 can be applied), for which the convergence results (i)-(ii)-(iii) of Proposition \ref{P:PenalizedBSDE} hold on $[0,T]$. Even though the maximal solution to \eqref{BSDE}-\eqref{BSDE_Constraint} on $[0,T]$ can a priori depends on $T$, this is not the case due to the convergences (i)-(ii)-(iii) on $[0,T]$, which call in the penalized BSDE \eqref{BSDEnk}, whose solution does not depend on $T$. As a consequence, we can past together all these maximal solutions on $[0,T]$, for any $T>0$, and we end up with a maximal solution to equation \eqref{BSDE}-\eqref{BSDE_Constraint} on $[0,\infty)$. Finally, the uniqueness of the maximal solution follows from Lemma \ref{L:Uniqueness}.
\ep

\subsubsection{Feynman-Kac formula for $V^{\beta,n}$}

We now derive, by means of the doubly indexed penalized BSDE \eqref{BSDEnk}, the Feynman-Kac formula for $V^{\beta,n}$ and study its viscosity properties. To this end, for any $\beta>0$, $n,k\in\N$, we define the function $V^{\beta,n,k}\colon\R^d\times\I_m\times U\rightarrow\R$ as follows
\begin{equation}
\label{Vbeta,n,k}
V^{\beta,n,k}(x,i,u) \ = \ Y_0^{\beta,n,k,x,i,u}, \qquad \forall\,(x,i,u)\in\R^d\times\I_m\times U.
\end{equation}
Then, $V^{\beta,n,k}$ is associated to the following elliptic integro-PDE:
\begin{align}
\label{HJB_penalizednk}
\beta V^{\beta,n,k}(x,i,u) - \Lc^{i,u}V^{\beta,n,k}(x,i,u) - \int_U \big[V^{\beta,n,k}(x,i,u') - V^{\beta,n,k}(x,i,u)\big] \vartheta_\mu(du') & \notag \\
- f(x,i,u) - n\sum_{j=1}^m \big[V^{\beta,n,k}(x,j,u) - V^{\beta,n,k}(x,i,u)
- c(x,i,j)\big]^+ & \\
+ k \int_U \big[V^{\beta,n,k}(x,i,u') - V^{\beta,n,k}(x,i,u)\big]^- \vartheta_\mu(du') & \ = \ 0, \notag
\end{align}
for any $(x,i,u)\in\R^d\times\I_m\times U$. More precisely, we have the following result.

\begin{Lemma}
Let Assumptions {\bf (H1)} and {\bf (H2)} hold. Then, for any $\beta>0$, $n,k\in\N$, the function $V^{\beta,n,k}$ defined in \eqref{Vbeta,n,k} is a continuous viscosity solution to \eqref{HJB_penalizednk} satisfying
\begin{equation}
\label{EstimateVbeta,n,k}
\big|V^{\beta,n,k}(x,i,u)\big| \ \leq \ \frac{C_{b,\sigma,f}}{\beta} \big(1 + |x|\big), \qquad \forall\,(x,i,u)\in\R^d\times\I_m\times U,
\end{equation}
where $C_{b,\sigma,f}$ is the same constant as in estimate \eqref{EstimateYbeta,n,k}.
\end{Lemma}
\textbf{Proof.}
The proof is standard and can be done along the same lines as in the proof of Proposition 3.2 in \cite{CFP14}.
\ep

\vspace{3mm}

We can finally state the main result of this appendix, namely the Feynman-Kac formula for $V^{\beta,n}$, which proves statement (ii) of Proposition  \ref{T:Feynman-Kac}.

\begin{Proposition}
\label{T:Feynman-Kac_App}
Let Assumptions {\bf (H1)}, {\bf (H2)}, and {\bf (H$U$)} hold. Then, for any $\beta>0$, $n\in\N$, $(x,i)\in\R^d\times\I_m$, $u\in\mathring U$, the random variable $Y_0^{\beta,n,x,i,u}$ is equal $\P$-a.s. to a constant independent of $u\in\mathring U$. Moreover, the function $V^{\beta,n}\colon\R^d\times\I_m\rightarrow\R$ given by
\[
V^{\beta,n}(x,i) \ = \ Y_0^{\beta,n,x,i,u}, \qquad \forall\,(x,i)\in\R^d\times\I_m,
\]
for any $u\in\mathring U$, is a continuous viscosity solution to system \eqref{HJB_penalized} and satisfies
\begin{equation}
\label{LinGrowthCondVbetan}
\sup_{(x,i)\in\R^d\times\I_m}\frac{|V^{\beta,n}(x,i)|}{1 + |x|} \ < \ \infty.
\end{equation}
\end{Proposition}
\textbf{Proof.}
The proof can be done as in \cite{CFP14}, Theorem 3.1. Here, we recall the main steps. Firstly, for any $T>0$, $\beta>0$, $n,k\in\N$, we consider the following parabolic integro-PDE on $[0,T]\times\R^m\times\I_m\times U$ in the unknown $w\colon[0,T]\times\R^d\times\I_m\times U\rightarrow\R$:
\begin{align}
\label{HJB_penalizednk_T}
\beta w - \frac{\partial w}{\partial t} - \Lc^{i,u}w - \int_U \big[w(t,x,i,u') - w(t,x,i,u)\big] \vartheta_\mu(du') & \notag \\
- f(x,i,u) - n\sum_{j=1}^m \big[w(t,x,j,u) - w(t,x,i,u)
- c(x,i,j)\big]^+ & \\
+ k \int_U \big[w(t,x,i,u') - w(t,x,i,u)\big]^- \vartheta_\mu(du') & \ = \ 0, \notag
\end{align}
with terminal condition $w(T,x,i,u)=V^{\beta,n,k}(x,i,u)$, $\forall\,(x,i,u)\in\R^d\times\I_m\times U$. Since the coefficients of system \eqref{HJB_penalizednk_T} are constant with respect to time $t$, we see that $V^{\beta,n,k}$ solves \eqref{HJB_penalizednk_T} in the viscosity sense. Then, we can apply Theorem 3.1 in \cite{KP12} to conclude that the limit $\lim_{k\rightarrow\infty}V^{\beta,n,k}(x,i,u)=\lim_{k\rightarrow\infty}Y_0^{\beta,n,k,x,i,u}=Y_0^{\beta,n,x,i,u}$ does not depend on $u$ in the interior $\mathring U$ of $U$. Notice that Theorem 3.1 in \cite{KP12} applies to equations with ``$\sup_{u\in U}$'' instead of ``$\inf_{u\in U}$'' as in \eqref{HJB_penalized}; however, simply observe that if $V^{\beta,n}$ is a viscosity solution to system \eqref{HJB_penalized}, then $-V^{\beta,n}$ solves a system with ``$\sup_{u\in U}$'' in place of ``$\inf_{u\in U}$'', for which we can use the results of Theorem 3.1 in \cite{KP12}.

The continuity of the function $V^{\beta,n}$ is a consequence of estimate \eqref{V^beta,n_Lipschitz}, which can be proved without relying on the (not yet proven) viscosity properties of $V^{\beta,n}$, but proceeding as in Corollary \ref{C:EstimatesV^beta,n}, where we used the dual representation formula \eqref{DualFormulaV^beta,n}. Moreover, from the monotone convergence of $(V^{\beta,n,k})_k$ towards $V^{\beta,n}$ as $k\rightarrow\infty$, together with estimate \eqref{EstimateVbeta,n,k}, we deduce the linear growth condition \eqref{LinGrowthCondVbetan} of $V^{\beta,n}$.

Finally, thanks again to Theorem 3.1 in \cite{KP12}, we have that, for any $T>0$, the function $V^{\beta,n}$ is a viscosity solution to the system of parabolic PDEs on $[0,T]\times\R^d\times\I_m$ in the unknown $w\colon[0,T]\times\R^d\times\I_m\rightarrow\R$:
\begin{equation}
\label{HJB_penalized_T}
\beta w - \frac{\partial w}{\partial t} - \inf_{u\in U}\big[\Lc^{i,u}w + f(x,i,u)\big] - n\sum_{j=1}^m \big[w(t,x,j) - w(t,x,i)
- c(x,i,j)\big]^+ \ = \ 0,
\end{equation}
with terminal condition $w(T,x,i)=V^{\beta,n}(x,i)$, $\forall\,(x,i)\in\R^d\times\I_m$. Since system \eqref{HJB_penalized_T} holds for every $T>0$, and $V^{\beta,n}$ is constant with respect to time, it follows that $V^{\beta,n}$ is indeed a viscosity solution to \eqref{HJB_penalized}.
\ep

\vspace{7mm}

\small
\bibliographystyle{plain}
\bibliography{biblio}

\begin{thebibliography}{10}

\bibitem{arilio98}
M~Arisawa and P.L. Lions.
\newblock On ergodic stochastic control.
\newblock {\em Comm. in Partial Differential Equations}, 23:2187--2217, 1998.

\bibitem{bardi_capuzzo-dolcetta}
M.~Bardi and I.~Capuzzo-Dolcetta.
\newblock {\em Optimal control and viscosity solutions of
  {H}amilton-{J}acobi-{B}ellman equations}.
\newblock Systems \& Control: Foundations \& Applications. Birkh\"auser Boston,
  Inc., Boston, MA, 1997.

\bibitem{BarlesSouganidisI}
G.~Barles and P.~E. Souganidis.
\newblock On the large time behavior of solutions of {H}amilton-{J}acobi
  equations.
\newblock {\em SIAM J. Math. Anal.}, 31(4):925--939 (electronic), 2000.

\bibitem{BarlesSouganidisII}
G.~Barles and P.~E. Souganidis.
\newblock Space-time periodic solutions and long-time behavior of solutions to
  quasi-linear parabolic equations.
\newblock {\em SIAM J. Math. Anal.}, 32(6):1311--1323 (electronic), 2001.

\bibitem{MR2336272}
Guy Barles and Espen~R. Jakobsen.
\newblock Error bounds for monotone approximation schemes for parabolic
  {H}amilton-{J}acobi-{B}ellman equations.
\newblock {\em Math. Comp.}, 76(260):1861--1893 (electronic), 2007.

\bibitem{BCP14}
E.~Bayraktar, A.~Cosso, and H.~Pham.
\newblock Robust feedback switching control: dynamic programming and viscosity
  solutions.
\newblock {\em SIAM J. Cont. Optim., to appear}, 2015.

\bibitem{MR2676760}
E.~Bayraktar and M.~Egami.
\newblock On the one-dimensional optimal switching problem.
\newblock {\em Math. Oper. Res.}, 35(1):140--159, 2010.

\bibitem{BS13}
E.~Bayraktar and M.~S{\^{\i}}rbu.
\newblock Stochastic {P}erron's method for {H}amilton-{J}acobi-{B}ellman
  equations.
\newblock {\em SIAM J. Control Optim.}, 51(6):4274--4294, 2013.

\bibitem{benfre92}
A.~Bensoussan and J.~Frehse.
\newblock On {B}ellman equations of ergodic control in $\mathbb{R}^n$.
\newblock {\em J. Reine Angew. Math}, 429:125--160, 1992.

\bibitem{CarLud10}
R.~Carmona and M.~Ludkovski.
\newblock Valuation of energy storage: an optimal switching approach.
\newblock {\em Quant. Finance}, 10(4):359--374, 2010.

\bibitem{CC16}
S.~Choukroun and A.~Cosso.
\newblock Backward {SDE} representation for stochastic control problems with
  nondominated controlled intensity.
\newblock {\em Ann. Appl. Probab.}, 26(2):1208--1259, 2016.

\bibitem{CCP15}
S.~Choukroun, A.~Cosso, and H.~Pham.
\newblock Reflected {BSDE}s with nonpositive jumps, and controller-and-stopper
  games.
\newblock {\em Stochastic Process. Appl.}, 125(2):597--633, 2015.

\bibitem{CFP14}
A.~Cosso, M.~Fuhrman, and H.~Pham.
\newblock Long time asymptotics for fully nonlinear {B}ellman equations: {A}
  backward {SDE} approach.
\newblock {\em Stochastic Process. Appl.}, 126(7):1932--1973, 2016.

\bibitem{CPH17}
A.~Cosso, H.~Pham, and H.~Xing.
\newblock {BSDE}s with diffusion constraint and viscous {H}amilton-{J}acobi
  equations with unbounded data.
\newblock {\em \textnormal{To appear on} Ann. Inst. Henri Poincar\'e Probab.
  Stat., \textnormal{preprint arXiv:1505.06868}}, 2017.

\bibitem{crandall_ishii_lions}
M.~G. Crandall, H.~Ishii, and P.-L. Lions.
\newblock User's guide to viscosity solutions of second order partial
  differential equations.
\newblock {\em Bull. Amer. Math. Soc. (N.S.)}, 27(1):1--67, 1992.

\bibitem{DucZer01}
K.~Duckworth and M.~Zervos.
\newblock A model for investment decisions with switching costs.
\newblock {\em Ann. Appl. Probab.}, 11(1):239--260, 2001.

\bibitem{ElHam09}
B.~El~Asri and S.~Hamad{\`e}ne.
\newblock The finite horizon optimal multi-modes switching problem: the
  viscosity solution approach.
\newblock {\em Appl. Math. Optim.}, 60(2):213--235, 2009.

\bibitem{EliKha14}
R.~Elie and I.~Kharroubi.
\newblock B{SDE} representations for optimal switching problems with controlled
  volatility.
\newblock {\em Stoch. Dyn.}, 14(3):1450003, 15, 2014.

\bibitem{fleming_souganidis}
W.~H. Fleming and P.~E. Souganidis.
\newblock On the existence of value functions of two-player, zero-sum
  stochastic differential games.
\newblock {\em Indiana Univ. Math. J.}, 38(2):293--314, 1989.

\bibitem{fuh-hu-tess}
M.~Fuhrman, Y.~Hu, and G.~Tessitore.
\newblock Ergodic {BSDE}s and optimal ergodic control in {B}anach spaces.
\newblock {\em SIAM J. Control Optim.}, 48:1542--1566, 2009.

\bibitem{HamJea07}
S.~Hamad{\`e}ne and M.~Jeanblanc.
\newblock On the starting and stopping problem: application in reversible
  investments.
\newblock {\em Math. Oper. Res.}, 32(1):182--192, 2007.

\bibitem{HamZha10}
S.~Hamad{\`e}ne and J.~Zhang.
\newblock Switching problem and related system of reflected backward {SDE}s.
\newblock {\em Stochastic Process. Appl.}, 120(4):403--426, 2010.

\bibitem{humadric14}
Y.~Hu, P.-Y. Madec, and A.~Richou.
\newblock A probabilistic approach to large time behavior of mild solutions of
  {HJB} equations in infinite dimension.
\newblock {\em SIAM J. Control Optim.}, 53(1):378--398, 2015.

\bibitem{HuTan10}
Y.~Hu and S.~Tang.
\newblock Multi-dimensional {BSDE} with oblique reflection and optimal
  switching.
\newblock {\em Probab. Theory Related Fields}, 147(1-2):89--121, 2010.

\bibitem{Hynd}
R.~Hynd.
\newblock The eigenvalue problem of singular ergodic control.
\newblock {\em Comm. Pure Appl. Math.}, 65(5):649--682, 2012.

\bibitem{ichihara12}
N.~Ichihara.
\newblock Large time asymptotic problems for optimal stochastic control with
  superlinear cost.
\newblock {\em Stochastic Process. Appl.}, 122(4):1248--1275, 2012.

\bibitem{ichiharaishii08}
N.~Ichihara and H.~Ishii.
\newblock Asymptotic solutions of {H}amilton-{J}acobi equations with
  semi-periodic hamiltonians.
\newblock {\em Comm. Partial Differential Equations}, 33(5):784--807, 2008.

\bibitem{ichsheu13}
N.~Ichihara and S.J. Sheu.
\newblock Large time behavior of solutions of {Hamilton-Jacobi-Bellman}
  equations with quadratic nonlinearity in gradients.
\newblock {\em SIAM J. Math. Anal.}, 45(1):279--306, 2013.

\bibitem{jacod_shiryaev}
J.~Jacod and A.~N. Shiryaev.
\newblock {\em Limit theorems for stochastic processes}, volume 288 of {\em
  Grundlehren der Mathematischen Wissenschaften [Fundamental Principles of
  Mathematical Sciences]}.
\newblock Springer-Verlag, Berlin, second edition, 2003.

\bibitem{khaetal10}
I.~Kharroubi, J.~Ma, H.~Pham, and J.~Zhang.
\newblock Backward {SDE}s with constrained jumps and quasi-variational
  inequalities.
\newblock {\em Ann. Probab.}, 38(2):794--840, 2010.

\bibitem{KP12}
I.~Kharroubi and H.~Pham.
\newblock Feynman-{K}ac representation for {H}amilton-{J}acobi-{B}ellman
  {IPDE}.
\newblock {\em Ann. Probab.}, 43(4):1823--1865, 2015.

\bibitem{lions_perthame86}
P.-L. Lions and B.~Perthame.
\newblock Quasivariational inequalities and ergodic impulse control.
\newblock {\em SIAM J. Control Optim.}, 24(4):604--615, 1986.

\bibitem{LyvPha07}
V.~Ly~Vath and H.~Pham.
\newblock Explicit solution to an optimal switching problem in the two-regime
  case.
\newblock {\em SIAM J. Control Optim.}, 46(2):395--426 (electronic), 2007.

\bibitem{menaldi_perthame_robin90}
J.-L. Menaldi, B.~Perthame, and M.~Robin.
\newblock Ergodic problem for optimal stochastic switching.
\newblock {\em J. Math. Anal. Appl.}, 147(2):512--530, 1990.

\bibitem{Nagai}
H.~Nagai.
\newblock Bellman equations of risk-sensitive control.
\newblock {\em SIAM J. Control Optim.}, 34(1):74--101, 1996.

\bibitem{ngopha}
M.~Ngo and H.~Pham.
\newblock Optimal switching for the pairs trading rule: A viscosity solutions
  approach.
\newblock {\em J. Math. Anal. Appl.}, 441(1):403--425, 2016.

\bibitem{pham09}
H.~Pham.
\newblock {\em Continuous-time stochastic control and optimization with
  financial applications}, volume~61 of {\em Stochastic Modelling and Applied
  Probability}.
\newblock Springer-Verlag, Berlin, 2009.

\bibitem{pham_zhang}
T.~Pham and J.~Zhang.
\newblock Two person zero-sum game in weak formulation and path dependent
  {B}ellman-{I}saacs equation.
\newblock {\em SIAM J. Control Optim.}, 52(4):2090--2121, 2014.

\bibitem{quenez_sulem}
M.-C. Quenez and A.~Sulem.
\newblock B{SDE}s with jumps, optimization and applications to dynamic risk
  measures.
\newblock {\em Stochastic Process. Appl.}, 123(8):3328--3357, 2013.

\bibitem{robxin13}
S.~Robertson and H.~Xing.
\newblock Large time behavior of solutions to semi-linear equations with
  quadratic growth in the gradient.
\newblock {\em SIAM J. Cont. Optim.}, 53(1):185--212, 2015.

\bibitem{sonzha09}
Q.~Song and Q.~Zhang.
\newblock An optimal pairs-trading rule.
\newblock {\em Automatica}, 49:3007--3014, 2013.

\bibitem{tang_li}
S.~J. Tang and X.~J. Li.
\newblock Necessary conditions for optimal control of stochastic systems with
  random jumps.
\newblock {\em SIAM J. Control Optim.}, 32(5):1447--1475, 1994.

\end{thebibliography}

\end{document}